\documentclass[11pt, oneside]{amsart}
\usepackage[text={5.58in,8.5in},centering,letterpaper,dvips]{geometry}
\usepackage[dvipsnames]{xcolor}
\usepackage{graphicx}
\usepackage{subcaption}
\usepackage{amsfonts}
\usepackage{epsf}
\usepackage{amssymb}
\usepackage{amsmath}
\usepackage{amscd}
\usepackage{tikz}
\usepackage{pdfpages}
\usepackage{fancyhdr}
\usepackage{setspace}
\usepackage{hyperref}
\usepackage[all]{xy}
\usetikzlibrary{matrix}
\usepackage{verbatim}
\usepackage{enumerate}

\theoremstyle{theorem}
\newtheorem{theorem}{Theorem}[section]
\newtheorem{proposition}[theorem]{Proposition}
\newtheorem{lemma}[theorem]{Lemma}
\newtheorem{question}[theorem]{Question}
\newtheorem{corollary}[theorem]{Corollary}

\newtheorem*{Powell}{Powell Conjecture}

\theoremstyle{definition}

\newtheorem{remark}[theorem]{Remark}

\newcommand{\Z}{\mathbb{Z}}
\newcommand{\C}{\mathcal{C}}
\newcommand{\N}{\mathbb{N}}
\newcommand{\Q}{\mathbb{Q}}
\newcommand{\R}{\mathbb{R}}

\newcommand{\D}{\mathcal D}
\newcommand{\Hs}{\mathcal R}

\newcommand{\wh}[1]{\widehat{#1}}

\newcommand{\A}{\alpha}
\newcommand{\n}{\beta}

\newcommand{\pd}{\partial}

\newcommand{\emp}{\emptyset}

\newcommand{\G}{\mathcal{G}}
\newcommand{\Pw}{\mathcal{P}}
\newcommand{\be}{\begin{enumerate}}
\newcommand{\ee}{\end{enumerate}}

\makeatletter
\def\@seccntformat#1{%
  \protect\textup{\protect\@secnumfont
    \ifnum\pdfstrcmp{subsection}{#1}=0 \bfseries\fi
    \csname the#1\endcsname
    \protect\@secnumpunct
  }%
}  
\makeatother

\makeatletter
\newtheorem*{rep@theorem}{\rep@title}
\newcommand{\newreptheorem}[2]{%
\newenvironment{rep#1}[1]{%
 \def\rep@title{#2 \ref{##1}}%
 \begin{rep@theorem}}%
 {\end{rep@theorem}}}
\makeatother

\newreptheorem{theorem}{Theorem}
\newreptheorem{lemma}{Lemma}
\newreptheorem{question}{Question}
\newreptheorem{corollary}{Corollary}
\newreptheorem{proposition}{Proposition}

\begin{document}

\rhead{\thepage}
\lhead{\author}
\thispagestyle{empty}


\raggedbottom
\pagenumbering{arabic}
\setcounter{section}{0}


\title{The Powell Conjecture and reducing sphere complexes}

\author{Alexander Zupan}
\address{Department of Mathematics, University of Nebraska-Lincoln, Lincoln, NE 68588}
\email{zupan@unl.edu}
\urladdr{http://www.math.unl.edu/~azupan2}

\begin{abstract}
The Powell Conjecture offers a finite generating set for the genus $g$ Goeritz group, the group of automorphisms of $S^3$ that preserve a genus $g$ Heegaard surface $\Sigma_g$, generalizing a classical result of Goeritz in the case $g=2$.  We study the relationship between the Powell Conjecture and the reducing sphere complex $\Hs(\Sigma_g)$, the subcomplex of the curve complex $\C(\Sigma_g)$ spanned by the reducing curves for the Heegaard splitting.  We prove that the Powell Conjecture is true if and only if $\Hs(\Sigma_g)$ is connected.  Additionally, we show that reducing curves that meet in at most six points are connected by a path in $\Hs(\Sigma_g)$; however, we also demonstrate that even among reducing curves meeting in four points, the distance in $\Hs(\Sigma_g)$ between such curves can be arbitrarily large.  We conclude with a discussion of the geometry of $\Hs(\Sigma_g)$.
\end{abstract}

\maketitle

\section{Introduction}

Let $Y$ be a closed, orientable 3-manifold and let $\Sigma \subset Y$ be a Heegaard surface for $Y$.  The \emph{Goeritz group} $\G(Y,\Sigma)$ is the set of orientation-preserving homeomorphisms $\varphi: Y \rightarrow Y$, considered up to isotopy, such that $\varphi(\Sigma) = \Sigma$.  In the case that $Y = S^3$ and $\Sigma_g$ is the standard genus $g$ Heegaard surface for $S^3$, the group $\G(S^3,\Sigma_g)$, or simply $\G_g$, is classically known as \emph{the genus $g$ Goeritz group}~\cite{goeritz}.  The Powell Conjecture offers a generating set for $\G_g$ in the case that $g \geq 3$, extending work of Goeritz to characterize $\G_2$.
\begin{Powell}\cite{powell}
For every $g \geq 3$, the genus $g$ Goeritz group $\G_g$ is generated by the five elements $\varphi_{\omega}$, $\varphi_{\eta}$, $\varphi_{\eta_{12}}$, $\varphi_{\nu}$, and $\varphi_{\theta}$ shown in Figure~\ref{pg}.
\end{Powell}

In~\cite{powell}, Powell claimed a proof of the Powell Conjecture, but in 2003 Scharlemann discovered that Powell's proof contains a fatal error.  Recently, Freedman and Scharlemann established the Powell Conjecture in the case $g = 3$~\cite{ScharFreed}; however, the conjecture remains open for $g \geq 4$.  (In a different direction, Freedman and Scharlemann have also noticed that one of Powell's conjectured generators is a product of the others, and thus is redundant~\cite{marty}.)

The goal of this paper is to better understand the Powell Conjecture from the perspective of the \emph{curve complex}.  The curve complex $\C(\Sigma)$ of a surface $\Sigma$ is a well-known space with a variety of connections to low-dimensional topology.  (See Section~\ref{prelim} for definitions.)  Given a Heegaard splitting $S^3 = H_1 \cup_{\Sigma_g} H_2$, the \emph{disk complex} $\D(H_i)$, $i=1,2$ is the subcomplex of $\C(\Sigma_g)$ induced by those curves that bound compressing disks in $H_i$.  Lastly, for genus $g \geq 3$, the \emph{reducing sphere complex} $\Hs(\Sigma_g)$ is the subcomplex of $\C(\Sigma_g)$ spanned by those curves that bound disks in both handlebodies $H_1$ and $H_2$; hence, these curves are reducing curves for the splitting.

It unknown in general if the reducing sphere complex is connected.  We relate the reducing sphere complex to the Powell Conjecture by proving the following theorem:

\begin{theorem}\label{connected}
The Powell Conjecture is true if and only if the reducing sphere complex $\Hs(\Sigma_g)$ is connected for all $g$.
\end{theorem}

The proof is an argument by induction on the genus $g$.  Explicitly, we prove that for a given $g \geq 3$, the genus $k$ Powell Conjecture is true for all $k \leq g$ if and only if the complex $\Hs(\Sigma_k)$ is connected for all $k \leq g$.  Using Freedman and Scharlemann's recent proof that the Powell Conjecture is true for $g=3$~\cite{ScharFreed}, we obtain

\begin{corollary}\label{3con}
The reducing sphere complex $\Hs(\Sigma_3)$ is connected.
\end{corollary}

For genus $g \geq 2$, curves in $\C(\Sigma)$ are connected by an edge whenever they have disjoint representatives.  For $g =1$, this definition is modified so that curves are connected by an edge if they intersect once.  In analogy, the definition of the reducing sphere complex $\Hs(\Sigma_g)$ may be modified to genus $g=2$, in which minimally intersecting non-homotopic curves meet in four points~\cite{ScharThomp03}, determining the edges of $\Hs(\Sigma_2)$.  The complex $\Hs(\Sigma_2)$ was introduced by Scharlemann in~\cite{Schar04} in order to give a modern proof of Goeritz's original argument~\cite{goeritz}.  Conversely, an argument that Goeritz's Theorem implies that $\Hs(\Sigma_2)$ is connected appears as Proposition 2.6 in~\cite{ScharThomp03}.  In Theorem 2.7 of~\cite{Schar01}, Scharlemann proves that the Powell Conjecture implies connectedness for a complex with vertices corresponding to complete collections of reducing spheres, which is closely related to the reducing sphere complex, and he remarks on page 408 that the converse ought to be true as well.  The diligent reader will note that~\cite{Schar01} was published before the error in Powell's work was discovered, so that the contingency on the Powell Conjecture is not included in the statement of Theorem 2.7 of~\cite{Schar01}.

These structures have also been examined for other 3-manifolds.  In the case that $\Sigma$ is a genus two Heegaard surface for an arbitrary 3-manifold $Y$,  $\Hs(\Sigma)$, which is also called the \emph{Haken sphere complex} in the literature, has been studied and characterized by Cho, Koda, and Seo in~\cite{arccomplexes} and by Cho and Koda in~\cite{hakenspheres18}, in which they prove the surprising fact that for the genus two Heegaard splitting of many lens spaces, $\Hs(\Sigma)$ is not connected.  Most recently, Cho and Koda completed the classification of the Goeritz groups of all 3-manifolds admitting a genus two Heegaard splitting~\cite{chokoda19}.

For the other two main results of the paper, we analyze reducing curves that meet in six or fewer points.  For any two curves $c$ and $c'$ in a surface $\Sigma$, let $\iota(c,c')$ denote their geometric intersection number.  Connectivity of the curve complex and disk complex can be proved by inducting on this intersection number.  In the spirit of proving the minimal cases of such an argument for $\Hs(\Sigma_g)$, we show

\begin{theorem}\label{6pts}
If $c$ and $c'$ are reducing curves for the Heegaard splitting $S^3 = H_1 \cup_{\Sigma_g} H_2$ such that $\iota(c \cap c') \leq 6$, then $c$ and $c'$ are contained in the same connected component of $\Hs(\Sigma_g)$.
\end{theorem}

Each complex has a natural path metric; we denote the distance between two curves by $d_{\ast}(c,c')$.  In the case of the curve complex $\C(\Sigma)$, distance is bounded above by a function of intersection number,
\[ d_{\C}(c,c') \leq 2 \log_2(\iota(c,c')) + 2.\]
A thorough discussion of this inequality appears in~\cite{notes}.  Moreover, a similar inequality relates distance and intersection number in the disk complex (see Lemma 2.1 of~\cite{hamenstadt}, for instance).  In contrast, we prove that surprisingly this is not true in $\Hs(\Sigma_g)$.

\begin{theorem}\label{bigdist}
For $g \geq 3$ and any $n \in \N$, there exist reducing curves $c,c' \in \Hs(\Sigma_g)$ such that
\be
\item $\iota(c,c') = 4$,
\item $d_{\C}(c,c') = d_{\D(H_i)}(c,c') = 2$,
and
\item $d_{\Hs}(c,c') > n.$
\ee
\end{theorem}

As a corollary, we obtain

\begin{corollary}\label{qi}
For $g \geq 3$, none of natural inclusions of $\Hs(\Sigma_g)$ into $\C(\Sigma_g)$ or $\D(H_i)$ is a quasi-isometric embedding.
\end{corollary}

The paper is organized as follows:  In Section~\ref{prelim}, we introduce the Powell generators and a space we call the \emph{Powell complex} $\Pw(\Sigma_g)$ to serve as an intermediary between $\G_g$ and $\Hs(\Sigma_g)$.  In Section~\ref{connect}, we relate Powell equivalence classes of $\G_g$ to connected components of $\Pw(\Sigma_g)$, which we in turn relate to components of $\Hs(\Sigma_g)$.  In Section~\ref{4p}, we prove that any reducing curves that meet four times are connected by a path in $\Hs(\Sigma_g)$, and in Section~\ref{6p} we strengthen the argument to prove Theorem~\ref{6pts}.  Finally, in Section~\ref{large} we prove the final result, Theorem~\ref{bigdist}.

\textbf{Acknowledgements}.  We are grateful to Marty Scharlemann for bringing this problem to our attention and for a number of helpful conversations and suggestions, including several of the arguments in Section~\ref{4p}.  We also thank Saul Schleimer for helpful conversations and for making us aware of the result in Theorem~\ref{fardist}.  Finally, we thank Abby Thompson for sharing her insights about the problem.  The author is supported by NSF grant DMS-1664578.

\section{Preliminaries}\label{prelim}

All manifolds are assumed to be compact and orientable.  For a subspace $P$ of a manifold $M$, we let $N(P)$ (resp. $\overline{N}(P)$) denote an open (resp. closed) regular neighborhood of $P$ in $M$.  Let $\Sigma$ be a compact surface.  A \emph{curve} in $\Sigma$ is a free homotopy class of an essential simple closed loop in $\Sigma$.  The \emph{curve complex} $\C(\Sigma)$ is a simplicial complex whose vertices correspond to curves in $\Sigma$, and whose $k$-cells correspond to subsets $(c_0,\dots,c_k)$ of $k+1$ curves in $\Sigma$ with pairwise disjoint representatives.  For two curves $c$ and $c'$ in $\Sigma$, we let $\iota(c,c')$ denote the geometric intersection number, the minimum number of intersections among representatives of $c$ and $c'$.  It is well-known that if $\Sigma$ admits a hyperbolic metric, then each curve has a unique geodesic representative, and pairs of these representatives realize geometric intersection number.  Suppose that $c$ and $c'$ are disjoint curves in $\Sigma$ and $e$ is an arc with such that $e \cap c$ is one of its endpoints and $e \cap c'$ is the other endpoint.  Then $\overline{N}(c \cup e \cup c')$ is a an embedded pair of pants in $\Sigma$ whose boundary is is the disjoint union of $c$, $c'$, and a third curve $c''$.  We say that $c''$ is the result of \emph{banding $c$ and $c'$ along $e$}.

Suppose now that $\Sigma \subset Y$ is a Heegaard surface for $Y$, so that $Y = H_1 \cup_{\Sigma} H_2$ for handlebodies $H_1$ and $H_2$.  The \emph{disk complex} of $H_i$, denoted $\D(H_i)$, is the full subcomplex of $\C(\Sigma)$ spanned by the curves in $\Sigma$ that bound compressing disks in $H_i$.  Finally, the \emph{reducing sphere complex} of $\Sigma$, denoted $\Hs(\Sigma)$, is defined to be the full subcomplex of $\C(\Sigma)$ spanned by curves that bound compressing disks in \emph{both} $H_1$ in $H_2$; in other words, $\Hs(\Sigma) = \D(H_1) \cap \D(H_2)$.  Observe that every curve in $\Hs(\Sigma)$ is the intersection of a reducing sphere for the splitting with the splitting surface $\Sigma$.  Abusing notation and terminology, we will often refer to curves and vertices interchangeably; if we say that a curve is a vertex, we mean that the curve corresponds to that vertex in the relevant complex.

The vertex set of each connected component of any of the above complexes is naturally a metric space using the path metric; the distance between two vertices is smallest number of edges in a path connecting them.  It is known that $\C(\Sigma)$ and $\D(H_i)$ are connected; however, it is an open problem whether $\Hs(\Sigma)$ is connected.  To bypass this issue, we define an extended metric on the entire complex by letting the distance between vertices in distinct components be $\infty$.  We denote the distances in $\C(\Sigma)$, $\D(H_i)$, and $\Hs(\Sigma)$ by $d_{\C(\Sigma)}$, $d_{\D(H_i)}$, and $d_{\Hs(\Sigma)}$, respectively.  When the surface $\Sigma$ is unambiguous, we omit it from this notation.

Although an element $\varphi \in \G_g$ is an automorphism of $S^3$, we will typically be interested in its restriction to $\Sigma_g$; thus, for ease of notation, we will use $\varphi$ in place of $\varphi|_{\Sigma_g}$.

\subsection{The Powell generators and the Powell complex} 
\begin{figure}[h!]
	\centering
	\includegraphics[width=.6\textwidth]{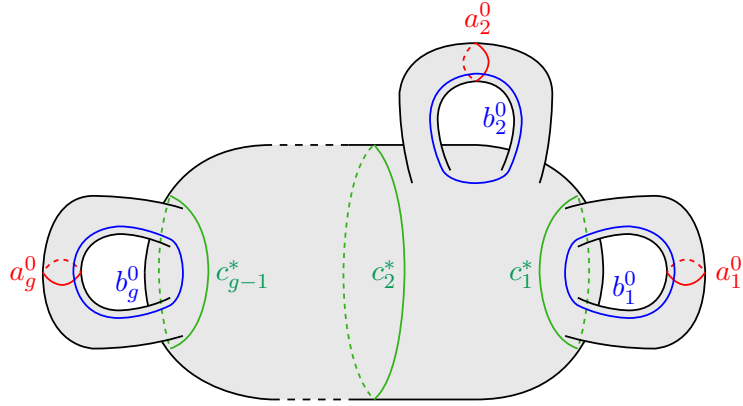}
\put (-268,50) {\textcolor{red}{$a_g^0$}}
\put (-228,46) {\textcolor{blue}{$b_g^0$}}
\put (-190,50) {\textcolor{ForestGreen}{$c_{g-1}^*$}}
\put (-131,50) {\textcolor{ForestGreen}{$c_2^*$}}
\put (-79,50) {\textcolor{ForestGreen}{$c_1^*$}}
\put (-1,50) {\textcolor{red}{$a_1^0$}}
\put (-40,44) {\textcolor{blue}{$b_1^0$}}
\put (-96,146) {\textcolor{red}{$a_2^0$}}
\put (-89,108) {\textcolor{blue}{$b_2^0$}}

	\caption{Curves in the standard diagram $v_0$}

	\label{standarddiagram}
\end{figure}

For the remainder of the paper, we suppose $S^3 = H_1 \cup_{\Sigma_g} H_2$ is the standard genus $g$ Heegaard splitting.  In~\cite{powell}, Powell proposes a set of generators for the genus $g$ Goeritz group $\G_g$.  The generators are defined relative to a fixed collection of $g$ genus one summands of the surface $\Sigma_g$.  Each summand contains a curve $a_i^0$ and bounds a disk in $H_1$ and a curve $b_i^0$ that bounds a disk in $H_2$, where $\iota(a_i^0,b_i^0) = 1$.  Let $\A^0 = \{a_1^0,\dots,a_g^0\}$, let $\n^0 = \{b_1^0,\dots,b_g^0\}$, and let $v_0 = (\A^0,\n^0)$.  The pair $v_0$ is an example of a \emph{standard diagram}, defined below.  We also fix a collection of reducing curves $c_1^*,\dots,c_{g-1}^*$, where $c_i^*$ separates the curves $\{a_1^0,b_1^0,\dots,a_i^0,b_i^0\}$ from the curves $\{a_{i+1}^0,b_{i+1}^0\,\dots,a_g^0,b_g^0\}$.  See Figure~\ref{standarddiagram}.\vspace{-.3cm}

\begin{figure}[h!]
\begin{subfigure}{.33\textwidth}
  \centering
  \includegraphics[width=.655\linewidth]{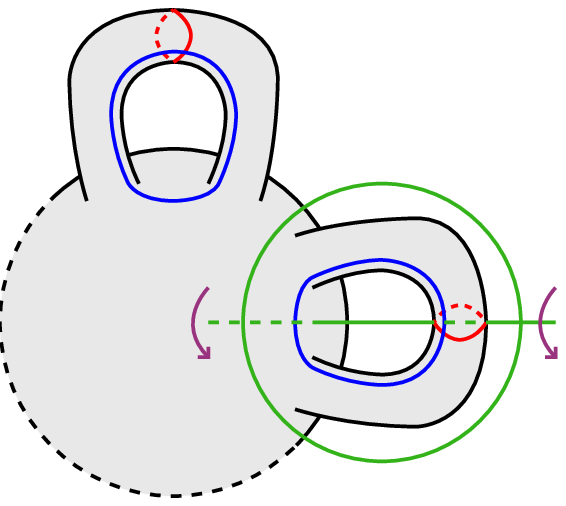}
  \caption{$\varphi_{\omega}$}
  \label{pgA}
\end{subfigure}%
\begin{subfigure}{.33\textwidth}
  \centering
  \includegraphics[width=.655\linewidth]{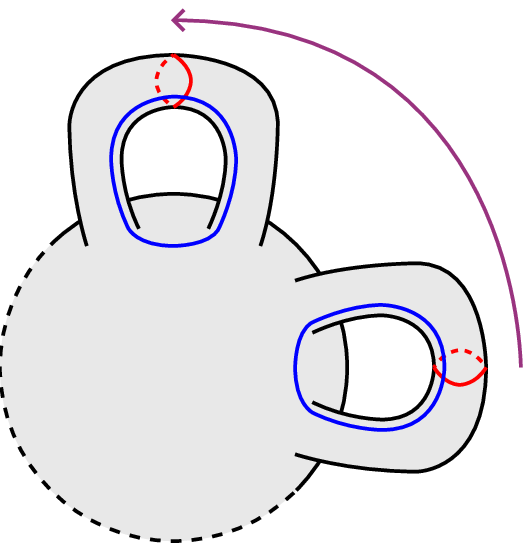}
  \caption{$\varphi_{\eta}$}
  \label{pgB}
\end{subfigure}
\begin{subfigure}{.33\textwidth}
  \centering
  \includegraphics[width=.655\linewidth]{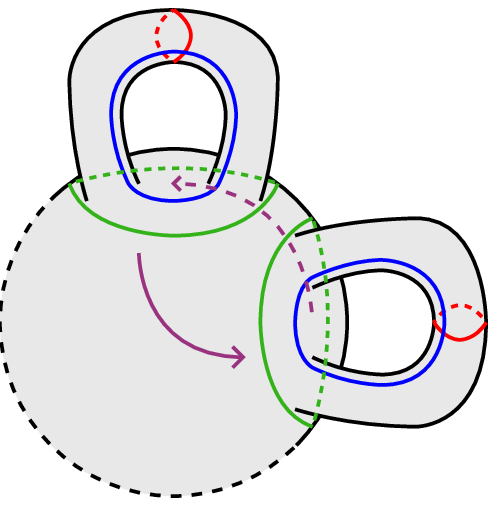}
  \caption{$\varphi_{\eta_{12}}$}
  \label{pgC}
\end{subfigure} \vspace{.3cm} \\

\begin{subfigure}{.24\textwidth}
  \centering
  \includegraphics[width=.9\linewidth]{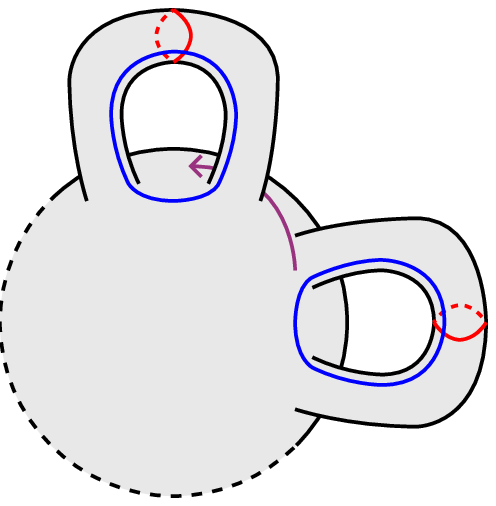}
  \caption{$\varphi_{\nu}[0]$}
  \label{pgD}
\end{subfigure}%
\begin{subfigure}{.24\textwidth}
  \centering
  \includegraphics[width=.9\linewidth]{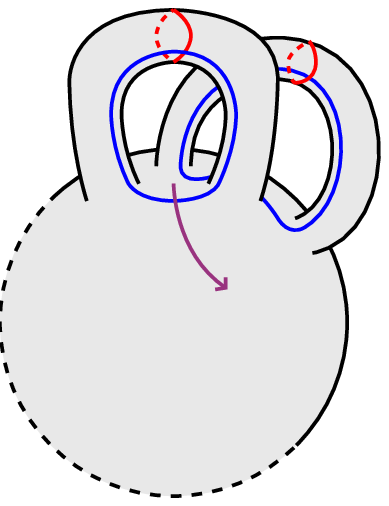}
  \caption{$\varphi_{\nu}[1/3]$}
  \label{pgE}
\end{subfigure}
\begin{subfigure}{.24\textwidth}
  \centering
  \includegraphics[width=.9\linewidth]{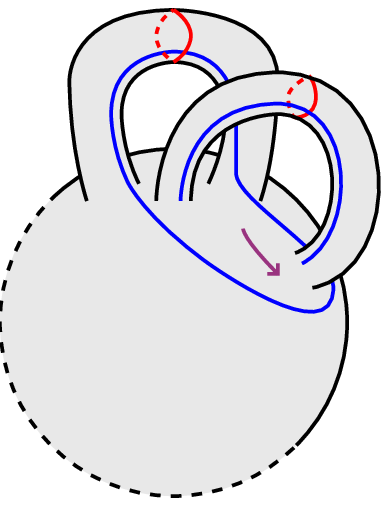}
  \caption{$\varphi_{\nu}[2/3]$}
  \label{pgF}
\end{subfigure}
\begin{subfigure}{.24\textwidth}
  \centering
  \includegraphics[width=.9\linewidth]{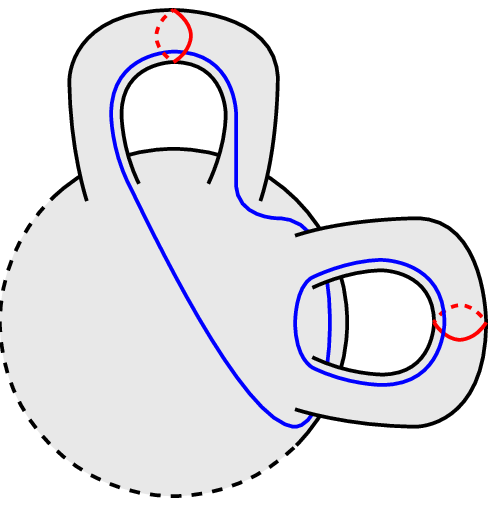}
  \caption{$\varphi_{\nu}[1]$}
  \label{pgG}
\end{subfigure}\vspace{.3cm} \\

\begin{subfigure}{.24\textwidth}
  \centering
  \includegraphics[width=.9\linewidth]{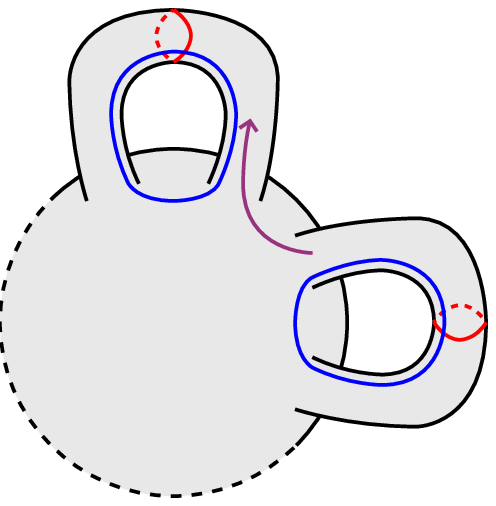}
  \caption{$\varphi_{\theta}[0]$}
  \label{pgH}
\end{subfigure}%
\begin{subfigure}{.24\textwidth}
  \centering
  \includegraphics[width=.9\linewidth]{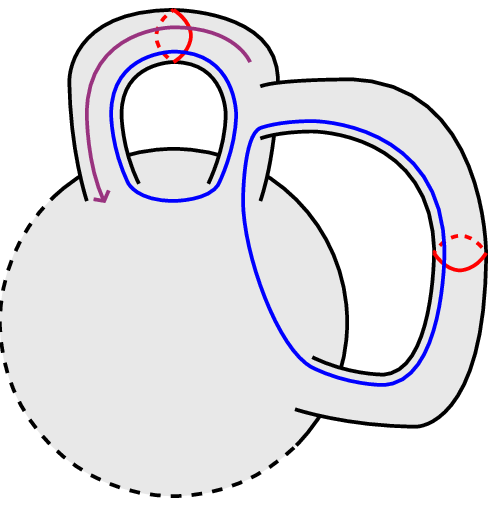}
  \caption{$\varphi_{\theta}[1/3]$}
  \label{pgI}
\end{subfigure}
\begin{subfigure}{.24\textwidth}
  \centering
  \includegraphics[width=.9\linewidth]{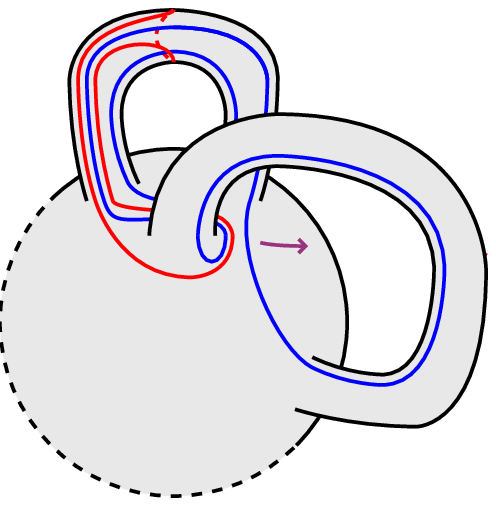}
  \caption{$\varphi_{\theta}[2/3]$}
  \label{pgJ}
\end{subfigure}
\begin{subfigure}{.24\textwidth}
  \centering
  \includegraphics[width=.9\linewidth]{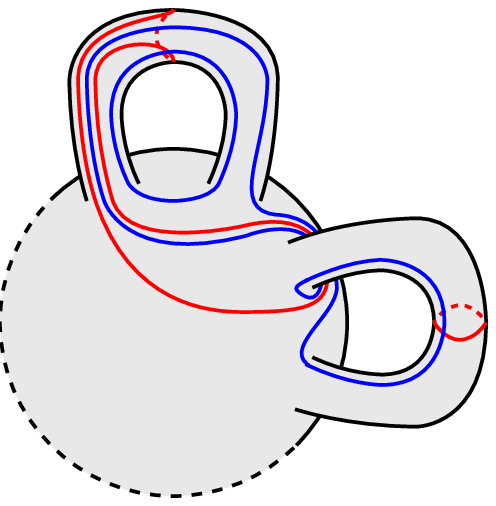}
  \caption{$\varphi_{\theta}[1]$}
  \label{pgK}
\end{subfigure}

\caption{Depictions of the Powell generators through excursions of $\Sigma_g$ in $S^3$.  In subfigures (D)-(K), we see four different snapshots of the two excursions yielding $\varphi_{\nu}$ and $\varphi_{\theta}$.  All pictured curves in $\A_0$ and $\n_0$ coincide with the curves $a_1^0$, $a_2^0$, $b_1^0$, and $b_2^0$, as labeled in Figure~\ref{standarddiagram}.}
\label{pg}
\end{figure}

Following~\cite{JM} and~\cite{ScharFreed}, let $\text{Diff}(S^3)$ denote the group of orientation-preserving diffeomorphisms of $S^3$, and let $\text{Diff}(S^3,\Sigma_g)$ denote the subgroup of $\text{Diff}(S^3)$ that maps $\Sigma_g$ to $\Sigma_g$.  Powell shows that $\G_g$ is a quotient of $\pi_1(\text{Diff}(S^3)/\text{Diff}(S^3,\Sigma_g))$ by a $\Z_2$ subgroup, and there is a natural projection map from this fundamental group to $\G_g$.  The motivation for this perspective is that elements of $\G_g$ can be viewed as end of an isotopy of $S^3$ that begins with the identity and returns $\Sigma_g$ to itself setwise.  Depictions of these generators are given in Figure~\ref{pg} (see also page 199 of~\cite{powell} or page 3 of~\cite{ScharFreed}), and their corresponding homeomorphisms $\varphi_{\omega}, \varphi_{\eta}, \varphi_{\eta_{12}}, \varphi_{\nu}, \varphi_{\theta} \in \G_g$ are called \emph{Powell generators}.

In each subfigures of in Figure~\ref{pg}, we keep track of the action of the homeomorphism $\varphi$ on the curves in $v_0$.  The generator $\varphi_{\omega}$ is an involution of the genus one summand of $\Sigma_g$ containing the curves $a_1^0$ and $b_1^0$.  The generator $\varphi_{\nu}$ permutes the curves in $v_0$ by sending $a_i^0$ to $a_{i+1}^0$ and $b_i^0$ to $b_{i+1}^0$.  The generator $\varphi_{\eta_{12}}$ swaps the curve $a_1^0$ with $a_2^0$ and the curve $b_1^0$ with $b_2^0$.  In each of the three cases, we see that setwise $\varphi_{\omega}(v_0) = \varphi_{\nu}(v_0) = \varphi_{\nu_{12}}(v_0) = v_0$.  The same statement does not hold for $\varphi_{\nu}$ or $\varphi_{\theta}$.

In order to keep track of the actions of the Powell generators on curves in $\Sigma_g$, we restrict our attention to collections of curves that behave similarly to the fixed standard diagram $v_0$.  As in~\cite{ScharFreed}, we say that two sets of pairwise curves $\A = \{a_1,\dots,a_g\}$ and $\n = \{b_1,\dots,b_g\}$ in $\Sigma_g$ are \emph{orthogonal} if $\iota(a_i,b_j) = \delta_{ij}$.    A \emph{standard diagram} in $\Sigma_g$ is a pair of orthogonal sets of curves $(\A,\n)$ such that
\begin{enumerate}
\item $a_i$ bounds a disk in $H_1$ for all $i$, and
\item $b_j$ bounds a disk in $H_2$ for all $j$.
\end{enumerate}

An element $\varphi$ of $\G_g$ is called a \emph{Powell move} if $\varphi$ can be expressed as a product of Powell generators, and (as in~\cite{ScharFreed}), two elements $\varphi,\varphi' \in \G_g$ are called \emph{Powell equivalent} if $\varphi'\circ \varphi^{-1}$ is a Powell move, in which case we write $\varphi \sim \varphi'$.  Equivalently, if $P_g$ is the subgroup of $\G_g$ generated by the Powell generators, the Powell equivalence classes correspond precisely with $P_g \backslash \G_g$, the set of right cosets of $P_g$ in $\G_g$.  Thus, the Powell Conjecture asserts that $P_g = \G_g$; equivalently, there exists only one Powell equivalence class.  The following lemma will be useful in our analysis.

\begin{lemma}\label{samesend}\cite[Lemma 1.7]{ScharFreed}
Suppose $\varphi,\varphi' \in \G_g$ and let $(\A,\n)$ be a standard diagram.  If $\varphi(\A) = \varphi'(\A) = \A_0$ (resp. $\varphi(\n) = \varphi'(\n) = \n_0$), then $\varphi \sim \varphi'$.
\end{lemma}

Instead of viewing Powell moves as excursions of a Heegaard surface, we shift our focus to the perspective of the curve complex and related structures.  First, we define a new complex, called the \emph{Powell complex}, and we present an equivalent formulation of the Powell Conjecture in this setting.  The vertices of $\Pw(\Sigma_g)$ are defined to be in one-to-one correspondence with standard diagrams in $\Sigma_g$.

There are straightforward modifications of a standard diagram to obtain a different standard diagram, and these moves constitute the edges of the $\Pw(\Sigma_g)$.  The edges connect vertices that have either $2g-1$ or $2g-2$ curves in common, with some additional constraints that require another definition:  Let $v = (\A,\n)$ be a vertex of $\Pw(\Sigma_g)$, with $\A = \{a_1,\dots,a_g\}$ and $\n = \{b_1,\dots,b_g\}$.  For each index $i$, let $c_i = \pd \overline{N}(a_i \cup b_i)$, so that $c_i$ bounds a disk in both $H_1$ and $H_2$.  We will call $c_i$ the reducing curve \emph{induced by} $a_i$ and $b_i$, noting that $c_i \in \Hs(\Sigma)$.

\begin{figure}[h!]
\begin{subfigure}{.5\textwidth}
  \centering
  \includegraphics[width=.8\linewidth]{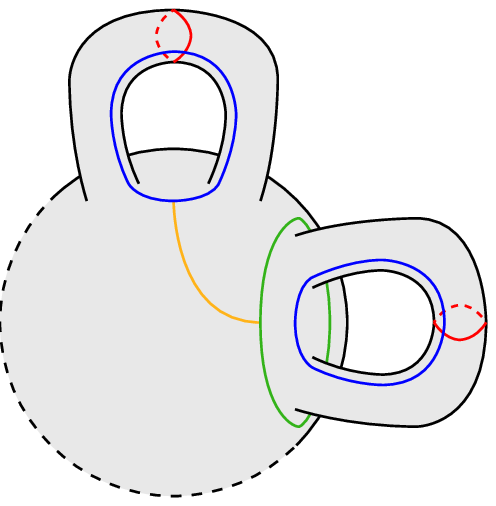}
\put (-120,65) {\textcolor{Dandelion}{$e$}}
\put (-101,40) {\textcolor{ForestGreen}{$c_i$}}
\put (-23,54) {\textcolor{red}{$a_i$}}
\put (-60,46) {\textcolor{blue}{$b_i$}}
\put (-123,154) {\textcolor{red}{$a_j$}}
\put (-114,114) {\textcolor{blue}{$b_j$}}
  \label{bob}
\end{subfigure}%
\begin{subfigure}{.5\textwidth}
  \centering
  \includegraphics[width=.8\linewidth]{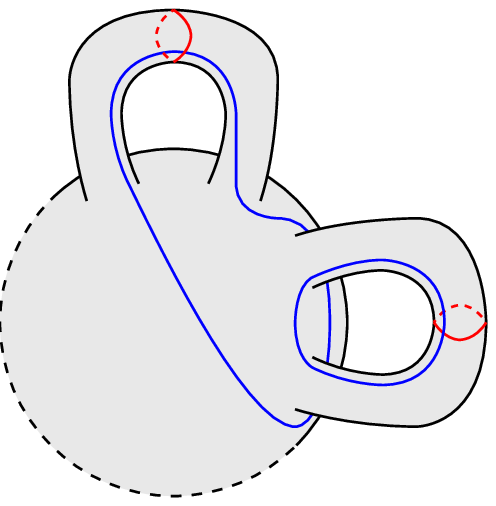}
\put (-23,54) {\textcolor{red}{$a_i$}}
\put (-60,46) {\textcolor{blue}{$b_i$}}
\put (-123,154) {\textcolor{red}{$a_j$}}
\put (-115,42) {\textcolor{blue}{$b'_j$}}
  \label{boo}
\end{subfigure}
  \caption{The result of a bubble move in which $b_j$ is replaced with the curve $b'_j$ obtained by banding $b_j$ to $c_i$}
\label{bub}
\end{figure}

For indices $i \neq j$, let $e$ be an arc in $\Sigma_g$ with endpoints in the curves $c_i$ and $a_j$ (resp. $b_j$) such that the interior of $e$ is disjoint from the curves in $v$, and let $a_j'$ (resp. $b_j'$) be the result of banding $a_j$ (resp. $b_j$) and $c_i$ along $e$.  Then the set $v'$ obtained from $v$ by replacing $a_j$ with $a_j'$ (resp. $b_j$ with $b_j'$) is again a standard diagram, and we say that $v$ and $v'$ are related by a \emph{bubble move}.  See Figure~\ref{bub}.

\begin{figure}[h!]
\begin{subfigure}{.5\textwidth}
  \centering
  \includegraphics[width=.8\linewidth]{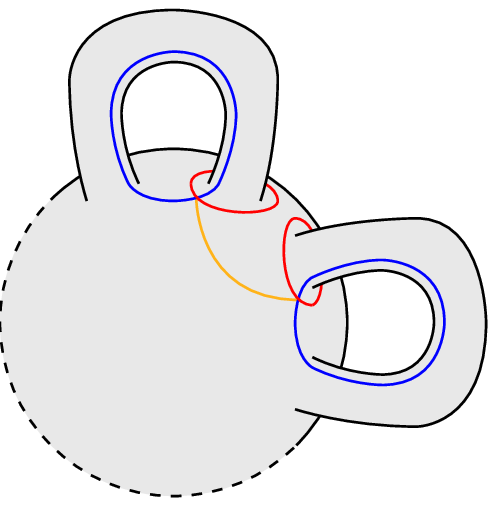}
\put (-113,70) {\textcolor{Dandelion}{$e$}}
\put (-94,68) {\textcolor{red}{$a_i$}}
\put (-60,46) {\textcolor{blue}{$b_i$}}
\put (-105,82) {\textcolor{red}{$a_j$}}
\put (-130,114) {\textcolor{blue}{$b_j$}}
  \label{eyb}
\end{subfigure}%
\begin{subfigure}{.5\textwidth}
  \centering
  \includegraphics[width=.8\linewidth]{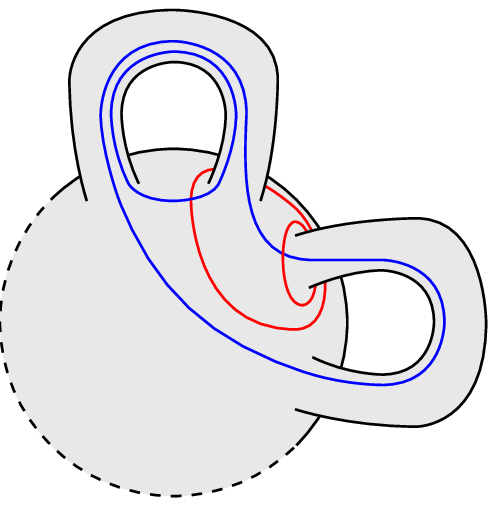}
\put (-94,65) {\textcolor{red}{$a_i$}}
\put (-109,82) {\textcolor{red}{$a_j'$}}
\put (-115,46) {\textcolor{blue}{$b'_i$}}
\put (-130,114) {\textcolor{blue}{$b_j$}}
  \label{eyo}
\end{subfigure}
  \caption{The result of an eyeglass move}
\label{eyeg}
\end{figure}

In another construction, let $i \neq j$, and consider an arc $e$ whose endpoints are the points $a_i \cap b_i$ and $a_j \cap b_j$ and such that the interior of $e$ is disjoint from the curves in $v$, and suppose further that the cyclic ordering of $a_i$, $b_i$, and $e$ obtained by traveling counterclockwise in a neighborhood of $a_i \cap b_i$ is opposite that of $a_j$, $b_j$, and $e$ in a neighborhood of $a_j \cap b_j$.  Let $a_i'$ be the curve obtained by banding $a_i$ to $a_j$ along $e$, and let $b_j'$ be the curve obtained by banding $b_j$ to $b_i$ along $e$.  Let $v'$ be the set of curves obtained from $v$ by replacing $a_j$ with $a_j'$ and $b_i$ with $b_i'$.  Then $v'$ is again a standard diagram and we say $v$ and $v'$ are related by an \emph{eyeglass move}.  The arc $e$ is called the \emph{bridge} of the eyeglass move, and the curves $a_j$ and $b_i$ are called the \emph{lenses}.  See Figure~\ref{eyeg}.

The edges of the Powell Complex $\Pw(\Sigma_g)$ are defined to correspond precisely to vertices related by bubble moves and eyeglass moves; as such, we will distinguish the two types of edges by calling them \emph{bubble edges} and \emph{eyeglass edges}, respectively.  There is a natural definition of higher-dimensional cells in the Powell Complex given by considering moves which commute, in the sense that the arcs defining their slides are disjoint and connect curves of distinct indices.  However, in this paper, we will only be concerned with the 1-skeleton of the Powell Complex.

Observe that the definition of bubble and eyeglass moves appear asymmetric; arguably, the edges in $\Pw(\Sigma_g)$ ought to be directed.  In fact, the moves are reversible, as demonstrated by the next lemma, and thus $\Pw(\Sigma_g)$ has unoriented edges.

\begin{lemma}\label{reverse}
Let $v$ and $v'$ be vertices in $\Pw(\Sigma_g)$.
\be
\item If there is a bubble move from $v$ to $v'$, then there is a bubble move from $v'$ to $v$.
\item If there is an eyeglass move from $v$ to $v'$, then there is an eyeglass move from $v'$ to $v$.
\ee
\end{lemma}

\begin{figure}[h!]
\begin{subfigure}{.5\textwidth}
  \centering
  \includegraphics[width=.8\linewidth]{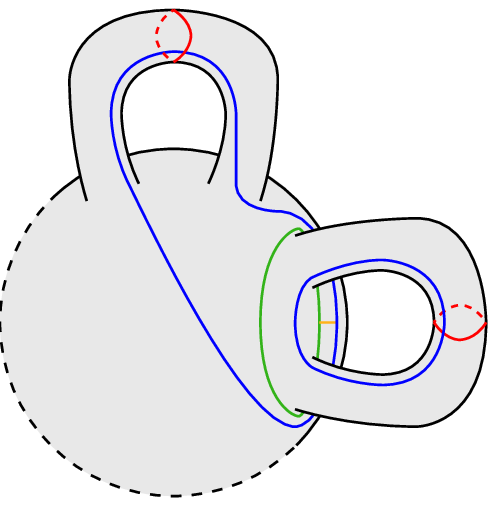}
  \put (-101,55) {\textcolor{ForestGreen}{$c_i$}}
\put (-81,53) {\textcolor{Dandelion}{$e'$}}
\put (-23,54) {\textcolor{red}{$a_i$}}
\put (-60,46) {\textcolor{blue}{$b_i$}}
\put (-123,154) {\textcolor{red}{$a_j$}}
\put (-115,42) {\textcolor{blue}{$b'_j$}}
  \label{eyb}
\end{subfigure}%
\begin{subfigure}{.5\textwidth}
  \centering
  \includegraphics[width=.8\linewidth]{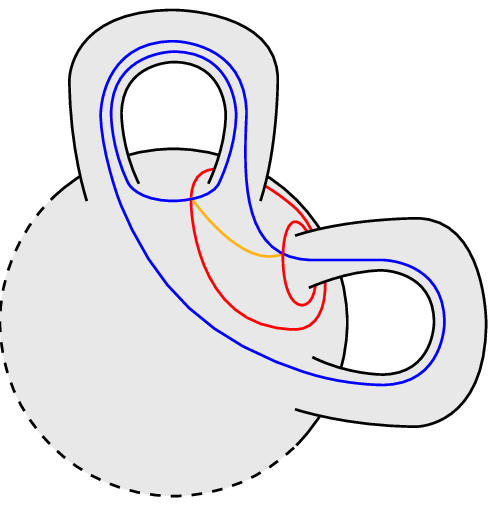}
\put (-94,65) {\textcolor{red}{$a_i$}}
\put (-122,82) {\textcolor{red}{$a_j'$}}
\put (-115,46) {\textcolor{blue}{$b'_i$}}
\put (-130,114) {\textcolor{blue}{$b_j$}}
\put (-105,86) {\textcolor{Dandelion}{$e'$}}
  \label{eyo}
\end{subfigure}
  \caption{Reversing a bubble move (left) and an eyeglass move (right).}
\label{revfig}
\end{figure}

\begin{proof}
For the first statement, the bubble move along arc $e'$ depicted at left in Figure~\ref{revfig} sends the right frame of Figure~\ref{bub}, showing $v'$, back to the left frame of Figure~\ref{bub}.  For the second statement, the eyeglass move with lenses $a_i$ and $b_j$ and arc $e'$ depicted at right in Figure~\ref{revfig} sends the right frame of Figure~\ref{eyeg}, showing $v'$, back to the left frame of Figure~\ref{eyeg}.
\end{proof}

Observe that the cyclic ordering of $a_i$, $b_i$, and $e$ near $a_i \cap b_i$ shown in the left panel of Figure~\ref{eyeg} is opposite that of $a_i$, $b_i'$, and $e'$ near $a_i \cap b_i'$, shown in the right panel of Figure~\ref{revfig}.  We call the former eyeglass move a \emph{right-handed eyeglass move} and the latter a \emph{left-handed eyeglass move}.  Lemma~\ref{reverse} implies that the reverse of a right-handed eyeglass move is left-handed, and vice versa.


Observe that if $\varphi$ is an automorphism of $\Sigma_g$, then for any standard diagram $v \in \Pw(\Sigma_g)$, the collection $\varphi(v)$ is also a standard diagram.  In the special case of our distinguished vertex $v_0$, we have

\begin{lemma}\label{genact}
The Powell generators $\varphi_{\omega}$, $\varphi_{\eta}$, and $\varphi_{\eta_{12}}$ fix $v_0$.  The vertices $v_0$ and $\varphi_{\nu}(v_0)$ are related by a bubble move.  The vertex $v_0$ is related to $\varphi_{\theta}(v_0)$ by a left-handed eyeglass move.
\end{lemma}

\begin{proof}
The first claim was discussed above and follows from Subfigures (A)-(C) of Figure~\ref{pg}.  The second claim follows from comparing Figure~\ref{pg} (G) to Figure~\ref{bub}.  For the third claim, the verification that $v_0$ and $\varphi_{\theta}(v_0)$ are related by an eyeglass move is shown in Figure~\ref{eyever} below, which is seen to be left-handed by comparing to the right panel of Figure~\ref{revfig}.
\end{proof}

\begin{figure}[h!]
\begin{subfigure}{.5\textwidth}
  \centering
  \includegraphics[width=.8\linewidth]{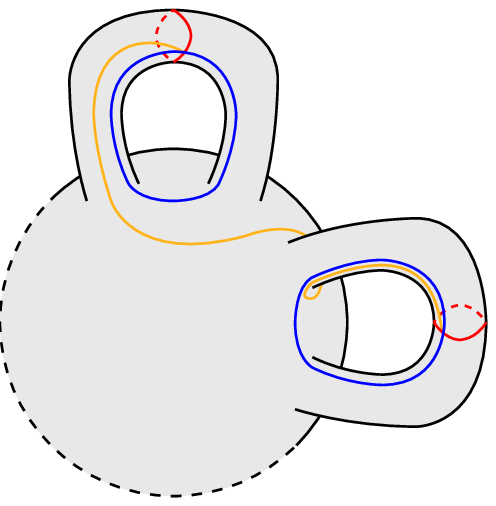}
\put (-113,69) {\textcolor{Dandelion}{$e^0$}}
\put (-23,54) {\textcolor{red}{$a_1^0$}}
\put (-60,46) {\textcolor{blue}{$b_1^0$}}
\put (-123,154) {\textcolor{red}{$a_2^0$}}
\put (-130,114) {\textcolor{blue}{$b_2^0$}}
  \label{eyb}
\end{subfigure}%
\begin{subfigure}{.5\textwidth}
  \centering
  \includegraphics[width=.8\linewidth]{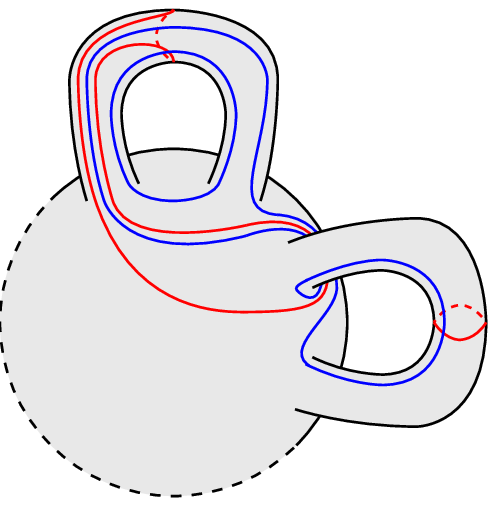}
\put (-23,54) {\textcolor{red}{$a_1^0$}}
\put (-60,46) {\textcolor{blue}{$b_1'$}}
\put (-123,154) {\textcolor{red}{$a_2'$}}
\put (-130,114) {\textcolor{blue}{$b_2^0$}}
  \label{eyo}
\end{subfigure}
  \caption{An eyeglass move representing the action of $\varphi_{\theta}$ on $v_0$.}
\label{eyever}
\end{figure}

Next, we note that for any $\varphi \in \G_g$, the vertices $v$ and $v'$ are connected by an edge in $\Pw(\Sigma_g)$ if and only if $\varphi(v)$ and $\varphi(v')$ are connected by an edge.  Therefore, $\varphi$ induces an automorphism (which we will also denote $\varphi$) of $\Pw(\Sigma_g)$.  We will let $\Pw_0$ denote the connected component of $\Pw(\Sigma_g)$ containing $v_0$.  In Section~\ref{connect}, we also use the action of $\varphi$ on $\Hs(\Sigma_g)$, which we record here as well.  We will let $\Hs_0$ denote the connected component of $\Hs(\Sigma_g)$ that contains the $g$ mutually disjoint reducing curves induced by $v_0$.

\begin{lemma}\label{pact}
If $\varphi$ is a Powell move, then $\varphi(\Pw_0) = \Pw_0$ and $\varphi(\Hs_0) = \Hs_0$.
\end{lemma}

\begin{proof}
The first claim follows immediately from Lemma~\ref{genact}.  For the second claim, let $c_i^0$ denote the reducing curve induced by $a_i^0$ and $b_i^0$.  Then we have $\varphi_{\omega}(c_i^0) = c_i^0$ for all $i$, $\varphi_{\eta}(c_i^0) = c_{i+1}^0$ for all $i$, and $\varphi_{\eta_{12}}(c_1^0) = c_2^0$, so each of these generators preserve $\Hs_0$.  In addition $\varphi_{\nu}(c_3^0) = \varphi_{\theta}(c_3^0) = c_3^0$, and we conclude that all generators and thus all Powell moves preserve $\Hs_0$.
\end{proof}

\begin{figure}[h!]
\begin{subfigure}{.33\textwidth}
  \centering
  \includegraphics[width=.8\linewidth]{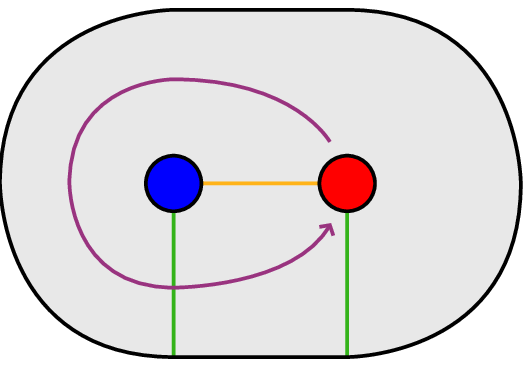}
  \caption{$\psi[0]$}
  \label{pgD}
\end{subfigure}%
\begin{subfigure}{.33\textwidth}
  \centering
  \includegraphics[width=.8\linewidth]{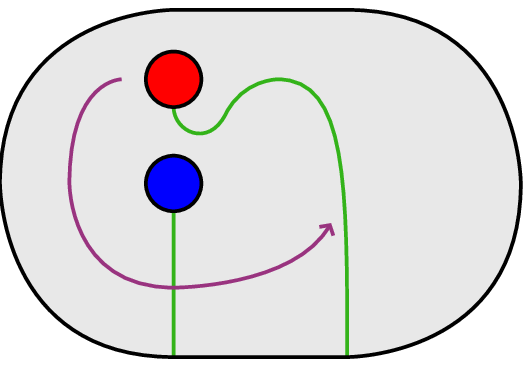}
  \caption{$\psi[1/3]$}
  \label{pgE}
\end{subfigure}
\begin{subfigure}{.33\textwidth}
  \centering
  \includegraphics[width=.8\linewidth]{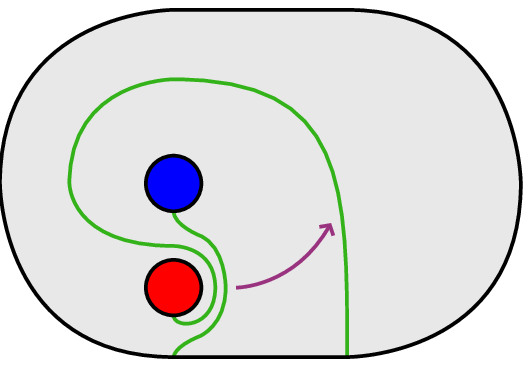}
  \caption{$\psi[2/3]$}
  \label{pgF}
\end{subfigure} \\

\begin{subfigure}{.33\textwidth}
  \centering
  \includegraphics[width=.8\linewidth]{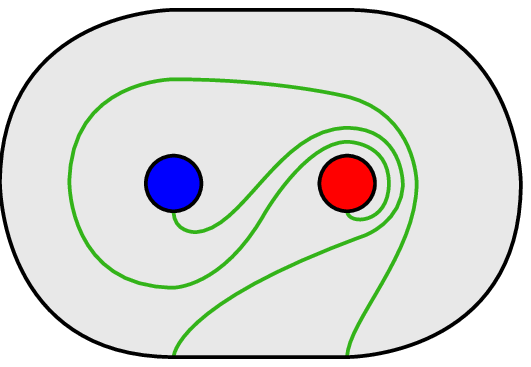}
  \caption{$\psi[1]$}
  \label{pgG}
\end{subfigure}
\begin{subfigure}{.33\textwidth}
  \centering
  \includegraphics[width=.8\linewidth]{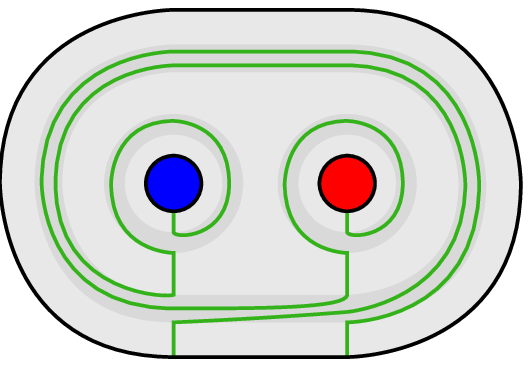}
  \caption{$\tau_a \circ \tau_b \circ \tau^{-1}_c$}
  \label{pgG}
\end{subfigure}
\caption{In (A)-(D), we see four different snapshots of an eyeglass twist.  Comparing (D) and (E), we verify that $\psi$ and $\tau_a \circ \tau_b \circ \tau_c^{-1}$ agree on the two arcs shown.}
\label{eyegdehn}
\end{figure}

The term \emph{eyeglass} relates to a particular type of element of $\G_g$ referred to by Freedman and Scharlemann as an \emph{eyeglass twist}~\cite{ScharFreed}:  Suppose that $a$ and $b$ are disjoint curves in $S$ that bound disks $D_a$ and $D_b$ in $H_1$ and $H_2$, respectively, and let $e$ be an embedded arc such that $a \cap e$ is one endpoint of $e$ and $b \cap e$ is the other.  Then $\pd \overline{N}(e \cup b)$ determines an arc $e'$ with endpoints in $a$, and an isotopy of $S^3$ that carries that disk $D_a$ counterclockwise around the arc $e'$ and back to its starting point yields an element of $\G_g$ called an \emph{eyeglass twist} with \emph{lenses} $a$ and $b$ and \emph{bridge} $e$ -- descriptively named because the union $a \cup e \cup b$ resembles a pair of eyeglasses.  See Figure~\ref{eyegdehn}.  Let $c$ be the curve resulting from banding $a$ and $b$ along $e$.  We call $c$ the \emph{boundary} of the eyeglasses.  Using this definition, we can see that the Powell generator $\varphi_{\theta}$ is an eyeglass twist with lenses $a^0_1$ and $b^0_2$ and bridge $e^0$ shown in Figure~\ref{eyever}.

A well-known homeomorphism of a surface is a \emph{Dehn twist}:  Let $c$ be a curve in a surface $\Sigma$, and parameterize a closed annular neighborhood $A$ of $c$ as $\{(e^{2\pi i s},t): s \in \R/\Z \text{ and } t \in [0,1]\}$.  The \emph{left-handed Dehn twist} $\tau_c: \Sigma \rightarrow \Sigma$ is defined to be the identity outside of $A$, and within $A$ it is given by $\tau_c(e^{2\pi i s},t) = (e^{2\pi i (s+t)},t)$.  For a comprehensive treatment, see~\cite{primer}.

\begin{lemma}\label{eyfact}
Suppose that $\psi \in \G_g$ is an eyeglass twist with lenses $a$ and $b$ and boundary $c$.  Then $\psi|_{\Sigma_g} = \tau_a \circ \tau_b \circ \tau^{-1}_c$.
\end{lemma}

\begin{proof}
Observe that $\psi$ is supported in a neighborhood of $a \cup e \cup b$, a pair of pants $P$.  Thus, it is generated by $\tau_a$, $\tau_b$, and $\tau_c$, and in addition, it is determined up to its action on a collection of arcs that cut $P$ into disks.  This action is depicted in Figure~\ref{eyegdehn}, from which we can deduce the desired statement.
\end{proof}

\begin{remark}
Using the factorization in Lemma~\ref{eyfact}, we note that $a$ bounds a disk in $H_1$ and $b \cup c$ bounds an annulus in $H_1$, so that the restriction of $\psi$ to $H_1$ is may be viewed as the product of a Dehn twist along the disk bounded by $a$ and an annulus twist along the annulus bounded by $b \cup c$.  A parallel argument can be used to understand the restriction of $\psi$ to $H_2$.
\end{remark}

\begin{remark}\label{lantern}
The interested reader may wish to compare Figure~\ref{eyegdehn} to Figure 5.3 from the well-known reference~\cite{primer}; in that setting, an eyeglass twist can be interpreted as a ``push map" of one boundary component of a pair of pants around the other.  Moreover, in the event that one of the disks, say $D_b$, can be split into two disks $D_{b_1}$ and $D_{b_2}$ disjoint from $D_a$, Freedman and Scharlemann demonstrate in Figure 8 of~\cite{ScharFreed} that the eyeglass twist with lenses $D_a$ and $D_b$ can be factored into the composition of an eyeglass twist with lenses $D_a$ and $D_{b_1}$ and another eyeglass twist $D_a$ and $D_{b_2}$.  This factorization is one justification for the famous \emph{lantern relation} in the mapping class group of a surface; this justification is discussed in Section 5.1.1 of~\cite{primer}
\end{remark} 

We can use Lemma~\ref{eyfact} to prove the following useful fact, connecting the terms ``lenses" and ``bridge" used for both eyeglass edges in $\Pw(\Sigma_g)$ and eyeglass twists in $\G_g$.

\begin{lemma}\label{eyenorm}
Suppose $v = (\A,\n)$ is related to $v' = (\A',\n')$ by a left-handed eyeglass move in $\Pw(T_g)$ with lenses $a_i \in \A$ and $b_j \in \n$ and bridge $e$.  Then there exists $\varphi \in \G_g$ such that the eyeglass twist $\psi$ with lenses $a_i$ and $b_j$ and bridge $e$ satisfies
\be
\item $\varphi(v) = v_0$,
\item $\varphi \circ \psi = \varphi_{\theta} \circ \varphi$, and
\item $\psi(v) = v'$.
\ee
\end{lemma}

\begin{proof}
Suppose $v =(\A,\n)$ and $v' = (\A',\n')$ are related by a left-handed eyeglass move with lenses $a_i$ and $b_j$ and bridge arc $e$.  Let $e^0$ be the arc depicted in Figure~\ref{eyever}.  By Lemma~\ref{genact}, $\varphi_{\theta}$ induces a left-handed eyeglass move on $v_0$, and so there is an orientation-preserving homeomorphism $\varphi:\Sigma_g \rightarrow \Sigma_g$ such that $\varphi(v) = v_0$, with the additional assumptions $\varphi(a_i) = a_1^0$, $\varphi(b_i) = b_1^0$, $\varphi(a_j) = a_2^0$, $\varphi(b_j) = b_2^0$, and $\varphi(e) = e^0$.  Since $\varphi$ extends over both $H_1$ and $H_2$, it determines an automorphism of $S^3$ and as such $\varphi \in \G_g$.

Let $\psi \in \G_g$ be an eyeglass twist with lenses $a_i$ and $b_j$ and bridge $e$, and let $c$ be the boundary of the eyeglasses.  By Lemma~\ref{eyfact}, we have $\psi = \tau_{a_i} \circ \tau_{b_j} \circ \tau^{-1}_c$.  Since the curves $a_j'$ and $b_i'$ in $v'$ are obtained by banding $a_i$ to $a_j$ along $e$ and $b_i$ to $b_j$ along $e$, respectively, it follows that the corresponding curves in $\varphi(v')$ are obtained by banding $\varphi(a_i)$ to $\varphi(a_j)$ and $\varphi(b_i)$ to $\varphi(b_j)$ along $\varphi(e)$; thus, we have $\varphi(v') = \varphi_{\theta}(v_0)$, as shown in Figure~\ref{eyever}.

As noted above, the generator $\varphi_{\theta}$ is an eyeglass twist with lenses $a_1^0$ and $b_2^0$ and bridge $e^0$; we let $c^0$ denote the boundary of the eyeglasses.  Again, by Lemma~\ref{eyfact}, we have $\varphi_{\theta} = \tau_{a_1^0} \circ \tau_{b_2^0} \circ \tau^{-1}_{c^0}$.  Therefore, by Fact 3.7 from~\cite{primer},
\[ \varphi \circ \psi = \varphi \circ (\tau_{a_i} \circ \tau_{b_j} \circ \tau_c^{-1}) = (\tau_{\varphi(a_i)} \circ \tau_{\varphi(b_j)} \circ \tau_{\varphi(c)}^{-1}) \circ \varphi = (\tau_{a_1^0} \circ \tau_{b_2^0} \circ \tau_{c^0}^{-1}) \circ \varphi = \varphi_{\theta} \circ \varphi.\]
It follows that $\varphi(v') = \varphi_{\theta}(v_0) = \varphi_{\theta}(\varphi(v)) = \varphi(\psi(v))$, and we conclude that $v' = \psi(v)$, completing the proof.
\end{proof}

Recall the reducing curves $c_1^*,\dots,c_{g-1}^*$ described above and shown in Figure~\ref{standarddiagram}.  The next lemma uses the factorization discussed in Remark~\ref{lantern} to show that a large family of eyeglass twists can be realized as Powell moves.

\begin{lemma}\cite[Lemma 3.4]{ScharFreed}\label{onebridge}
Suppose that $\psi$ is an eyeglass twist with bridge $e$ that intersects one of the reducing curves $c^*_i$ in a single point.  Then $\psi$ is a Powell move.
\end{lemma}

\section{Powell equivalence classes, $\Pw(\Sigma_g)$, and $\Hs(\Sigma_g)$}\label{connect}

In this section, we use the tools established thus far to relate Powell equivalence classes of elements of $\G_g$ to the connected components of $\Pw(\Sigma_g)$, which we in turn relate to the connected components of $\Hs(\Sigma_g)$.  

\begin{proposition}\label{pcomp}
The Powell equivalence classes of $\G_g$ are in one-to-one correspondence with the connected components of $\Pw(\Sigma_g)$.
\end{proposition}

\begin{proof}
We define a function $\Phi$ from the connected components of $\Pw(\Sigma_g)$ to the Powell equivalence classes of $\G_g$ and prove that $\Phi$ is a bijection.  Let $[v]$ represent the connected component of $\Pw(\Sigma_g)$ containing $v$, and let $[\varphi]$ represent the Powell equivalence class of $\G_g$ containing $\varphi$.  Recall that $\Pw_0$ is defined to be $[v_0]$.  For any vertex $v = (\A,\n) \in \Pw(\Sigma_g)$, let $\varphi$ be an automorphism of $\Sigma_g$ such that $\varphi(v) = v_0$.   Since $\varphi$ maps the Heegaard diagram $(\A,\n)$ to $(\A_0,\n_0)$, it follows that $\varphi$ extends to an automorphism of the pair $(S^3,\Sigma_g)$; hence $\varphi \in \G_g$.  Define $\Phi([v]) = [\varphi]$.  First, we show that $\Phi$ is well-defined.  If $\varphi'$ is another element of $\G_g$ such that $\varphi'(v) = v_0$, then $[\varphi] = [\varphi']$ by Lemma~\ref{samesend}.  Now, suppose that $v = (\A,\n)$ and $v' = (\A',\n')$ are connected by a bubble edge in $\Pw(\Sigma_g)$, and let $\varphi(v) = v_0$ and $\varphi'(v') = v_0$.  Since $v$ and $v'$ have $2g-1$ curves in common, either $\A = \A'$, in which case $\varphi(\A) = \varphi'(\A') = \A_0$, or $\n = \n'$, in which case $\varphi(\n) = \varphi'(\n') = \n_0$.  In either case, again invoking Lemma~\ref{samesend}, we have that $[\varphi] = [\varphi']$.

Next, suppose that $v = (\A,\n)$ and $v'= (\A',\n')$ are connected by an eyeglass edge in $\Pw(\Sigma_g)$.  If necessary, by Lemma~\ref{reverse} we can reverse the roles of $v$ and $v'$; thus, we may suppose without loss of generality that $v$ is related to $v'$ by a left-handed eyeglass move.  By Lemma~\ref{eyenorm}, there is an automorphism $\varphi$ of $\Sigma_g$ and an eyeglass twist $\psi$ such that $\varphi(v) = v_0$, $\varphi \circ \psi = \varphi_{\theta} \circ \varphi$, and $\psi(v) = v'$.  Rearranging yields $\varphi_{\theta} = \varphi \circ \psi \circ \varphi^{-1} = \varphi \circ (\varphi \circ \psi^{-1})^{-1}$; hence $\varphi \sim (\varphi \circ \psi^{-1})$.  Finally, note that $\psi(v) = (v')$ and so $\varphi(\psi^{-1}(v')) = \varphi(v) = v_0$.  We conclude that $\Phi([v']) = [\varphi \circ \psi^{-1}] = [\varphi] = \Phi([v])$.  Since every vertex in $[v]$ is connected to $v$ by a sequence of bubble and eyeglass edges, we have that $\Phi([v])$ is well-defined.

To see that $\Phi$ is surjective, let $\varphi \in \G_g$.  Then $\varphi^{-1}(v_0)$ is a vertex in $\Pw(\Sigma_g)$ and we have $\Phi([\varphi^{-1}(v_0)]) = [\varphi]$.  To complete the proof, suppose that $v,v' \in \Pw(\Sigma_g)$ and $\Phi([v]) = \Phi([v'])$.  Then there are elements $\varphi,\varphi' \in \G_g$ such that $\varphi(v) = v_0 = \varphi'(v')$, and $\varphi' \circ \varphi^{-1}$ is a Powell move.  We wish to show that $v$ is connected to $v'$ by a path in $\Pw(\Sigma_g)$.  By Lemma~\ref{pact}, we have that $v_0$ and $\varphi'(\varphi^{-1}(v_0)) = \varphi'(v)$ are contained in the same connected component, $\Pw_0$, of $\Pw(\Sigma_g)$.  Since $(\varphi')^{-1}$ induces an automorphism on $\Pw(\Sigma_g)$, it follows that $(\varphi')^{-1}(v_0) = v'$ and $(\varphi')^{-1}(\varphi'(v)) = v$ are also contained in the same connected component, $(\varphi')^{-1}(\Pw_0)$, of $\Pw(\Sigma_g)$.  We conclude that $[v] = [v']$, and $\Phi$ is a bijection.
\end{proof}

We remark that it is well-known that $\Pw(\Sigma_1)$ is a single vertex, and by Goeritz's classical theorem~\cite{goeritz} and Proposition~\ref{pcomp}, the complex $\Pw(\Sigma_2)$ is connected.

The next step in the argument is finding a relationship between components of $\Pw(\Sigma_g)$ and components of $\Hs(\Sigma_g)$; as above, if $c$ is a reducing curve, we let $[c]$ denote the component of $\Hs(\Sigma_g)$ containing $c$.  Given a vertex $v \in \Pw(\Sigma_g)$ we note that its collection of induced reducing curves $\{c_1,\dots,c_g\}$ are pairwise disjoint; hence we define a function $\Psi$ from connected components of $\Pw(\Sigma_g)$ to connected components of $\Hs(\Sigma_g)$ by the rule $\Psi([v]) = [c_i]$.  We split our analysis of $\Psi$ into two different propositions.

\begin{proposition}\label{psisurjective}
The function $\Psi$ is well-defined and surjective.
\end{proposition}

\begin{proof}
To see that $\Psi$ is well-defined, consider two vertices $v = (\A,\n)$ and $v' =(\A',\n')$ connected by a single edge in $\Pw(\Sigma_g)$.  Since $g \geq 3$, there is a pair of curves $(a_i,b_i) \in (\A \cap \A',\n \cap \n')$, with induced reducing curve $c_i$.  It follows that $\Psi([v]) = \Psi([v']) = [c_i]$, and since any two vertices in $[v]$ are connected by a sequence of edges, we see that $\Psi$ is well-defined.

To prove surjectivity, let $c$ be a reducing curve in $\Sigma_g$.  Then $c$ is the intersection of a reducing sphere $P$ with the Heegaard surface $\Sigma_g$, which can be reduced to smaller genus Heegaard surfaces $\Sigma_{g_1}$ and $\Sigma_{g_2}$ for $S^3$.  Each of these splittings has its own standard diagram $v_1$ and $v_2$.  Let $v = v_1 \cup v_2$.  Then $c$ is disjoint from any reducing curve induced by $v_1$ or $v_2$, and we have $\Psi([v]) = [c]$.
\end{proof}


In order to show injectivity, we need to strengthen our hypotheses and prove an additional lemma.

\begin{lemma}\label{injkey}
Suppose that the Powell complex $\Pw(\Sigma_k)$ is connected for all $k$ with $3 \leq k < g$.  Then any two vertices $v,v' \in \Pw(\Sigma_g)$ such that there exists a reducing curve $c$ with $v \cap c = v' \cap c = \emp$ satisfy $[v] = [v']$
\end{lemma}

\begin{proof}
Cutting $\Sigma_g$ along the reducing curve $c$ and capping the components with disks $D_1$ and $D_2$ naturally associates $\Sigma_g$ with the disjoint union of $\Sigma_{g_1}$ and $\Sigma_{g_2}$, where $g_1 + g_2 = g$.  Furthermore, letting $v_i = v \cap \Sigma_{g_i}$ and $v_i' = v' \cap \Sigma_{g_i}$, we have that $v_i,v_i' \in \Pw(\Sigma_{g_i})$.  By assumption, $\Pw(\Sigma_{g_i})$ is connected, so that there are paths from $v_i$ to $v_i'$.  Generically, the arcs yielding each of the bubble or eyeglass moves in these paths can be chosen to be disjoint from the caps $D_1$ and $D_2$.  If follows that $v$ is connected to a vertex $v''$ in $\Pw(\Sigma_g)$ such that $v'' \cap c = \emp$, $v''$ splits into $v''_1 \cup v''_2$, and the curves of $v''_i$ are isotopic to the curves of $v'_i$ in $\Sigma_{g_i}$.  

Note, however, that we do not necessarily know that curves in $v''_i$ are isotopic to $v'_i$ in $\Sigma_g$, since the isotopy in $\Sigma_{g_i}$ might pass a curve over the cap $D_i$.  Nevertheless, we may realize an isotopy of a curve $a''$ (or $b''$) in $v''_1$ over $D_1$ by banding $a''$ (or $b''$) to $c$ in $\Sigma_g$, which in turn is equivalent to banding $a''$ (or $b''$) to each of the reducing curves induced by $v''_2$; that is, a sequence of bubble moves.  Furthermore, we may realize an isotopy of a curve $a''$ (or $b''$) in $v''_2$ over $D_2$ by banding $a''$ (or $b''$) to $c$ in $\Sigma_g$, which in turn is equivalent to banding $a''$ (or $b''$) to each of the reducing curves induced by $v'_1$, another sequence of bubble moves.  We conclude that $v''$ is connected to $v_1' \cup v_2''$, which is in turn connected to $v'$ by a sequence of bubble edges in $\Pw(\Sigma_g)$, and thus $[v] = [v''] = [v']$.
\end{proof}

\begin{proposition}\label{psiinjective}
Suppose that the Powell complex $\Pw(\Sigma_k)$ is connected for all $k$ with $3 \leq k < g$.  Then $\Psi$ is injective.
\end{proposition}

\begin{proof}
Suppose that $v$ and $v'$ are vertices in $\Pw(\Sigma_k)$ such that $\Psi([v]) = \Psi([v'])$.  Fix a reducing curve $c$ induced by $v$ and a reducing curve $c'$ induced by $v'$.  Then $\Psi([v]) = [c] = [c'] = \Psi([v'])$, and there exists a path of reducing curves $c = c_0,c_1,\dots,c_n = c'$ in $\Hs(\Sigma_g)$.  Let $v_0 = v$, and $v_n = v'$, and for each $i$ such that $0 < i < n$, choose a vertex $v_i \in \Pw(\Sigma_g)$ such that $v_i \cap c_i = \emp$.  Additionally, for each index $i$ such that $0 \leq i < n$, the reducing curves $c_i$ and $c_{i+1}$ are disjoint; by splitting $\Sigma_g$ into three summands and choosing standard diagrams in each summand, we can choose a vertex $w_i \in \Pw(\Sigma_g)$ such that $w_i \cap c_i  = w_i \cap c_{i+1} = \emp$.

Note that for all $i$ with $0 \leq i < n$, we have $v_i \cap c_i = w_i \cap c_i = \emp$; hence, by Lemma~\ref{injkey}, $[v_i] = [w_i]$.  Likewise, $w_i \cap c_{i+1} = v_{i+1} \cap c_{i+1} = \emp$, so again by Lemma~\ref{injkey}, $[w_i]=[v_{i+1}]$ are in the same connected component of $\Pw(\Sigma_g)$.  We conclude that $[v_i] = [v_{i+1}]$ for all $i$, and in particular $[v] = [v_0] = [v_n] = [v']$, completing the proof.
\end{proof}

Finally, we combine these propositions to prove the first main theorem.

\begin{proof}[Proof of Theorem~\ref{connected}]
Fix $g \geq 3$.  We prove that the following three statements are equivalent:
\begin{enumerate}
\item The genus $k$ Goeritz Conjecture is true for all $k \leq g$.
\item The Powell complex $\Pw(\Sigma_k)$ is connected for all $k \leq g$.
\item The reducing sphere complex $\Hs(\Sigma_k)$ is connected for all $k \leq g$.
\end{enumerate}

First, Proposition~\ref{pcomp} implies that (1) and (2) are equivalent.  If $\Pw(\Sigma_g)$ has one connected component, then Proposition~\ref{psisurjective} implies that $\Hs(\Sigma_g)$ is connected as well.  To see that (3) implies (2), suppose that $\Hs(\Sigma_k)$ is connected for all $k \leq g$, and suppose by way of induction that (3) implies (2) for genera smaller than $g$, so that $\Pw(\Sigma_k)$ is connected for all $k < g$.  Propositions~\ref{psisurjective} and~\ref{psiinjective} imply that the function $\Psi$ is a bijection, so that $\Pw(\Sigma_g)$ and $\Hs(\Sigma_g)$ have the same number of connected components -- namely, one.
\end{proof}

\section{Reducing curves meeting in at most four points}\label{4p}

As discussed in the introduction, the well-known proofs that the curve complex and disk complex are connected induct on the intersection number of two curves to find a path between them.  In the section, we prove that if $c$ and $c'$ are reducing curves such that $\iota(c,c') \leq 4$, then $c$ and $c'$ are connected by a path in $\Hs(\Sigma_g)$, which we use in turn to prove the full generality of Theorem~\ref{6pts} in Section~\ref{6p}.

To this end, for the remainder of this section, we suppose $\iota(c,c') \leq 4$, we let $P$ and $P'$ denote the reducing spheres such that $c = P \cap \Sigma_g$ and $c' = P' \cap \Sigma_g$, and we assume that $P$ and $P'$ have been isotoped (fixing $c$ and $c'$) to intersect minimally in $S^3$.  In particular, this implies that every component of $P \cap P'$ meets the Heegaard surface $\Sigma_g$.  Suppose that $P = D_1 \cup D_2$ and $P' = D_1' \cup D_2'$, where $D_1$ and $D_1'$ are compressing disks in $H_1$, and $D_2$ and $D_2'$ are compressing disks in $H_2$.  Since each component of $P \cap P'$ meets $\Sigma_g$ and $\iota(c,c') \leq 4$, the intersection $P \cap P'$ contains either one or two curves.  

Suppose that $D$ and $D'$ are any two disks in one of the handlebodies $H_i$, and $D$ and $D'$ have been isotoped to intersect minimally, so that $D \cap D'$ is a collection of arcs.  If $\delta \subset D \cap D'$ is an arc that is outermost in $D$, then $\delta$ cuts out a subdisk $\Delta$ of $D$ whose interior is disjoint from $D'$.  In $D'$, the arc $\delta$ cuts $D'$ into two subdisks, $\Delta^+$ and $\Delta^-$, and we can obtain two new disks $D^+ = \Delta^+ \cup \Delta$ and $D^- = \Delta^- \cup \Delta$, pushing $D^{\pm}$ off of the disk $D$ along a collar neighborhood of the subdisk $\Delta$.  We say that $D^+$ and $D^-$ are obtained from $D'$ by \emph{surgery} along $\Delta$.  If either boundary component $\pd D^+$ or $\pd D^-$ is inessential, there exists an isotopy reducing $|D \cap D'|$; thus, it follows that both $D^+$ and $D^-$ are compressing disks for $H_i$.  Note that $D^{\pm} \cap D' = \emp$ and $|\pd D^+ \cap \pd D| + |\pd D^- \cap \pd D| < \iota(\pd D',\pd D)$; however, that it may be the case that neither disk $D^{\pm}$ intersects $D$ minimally.


We break the work in this section into three short lemmas and a more elaborate proposition.

\begin{lemma}\label{2pts}
If $\iota(c,c')= 2$, then $d_{\Hs}(c,c') = 2$.
\end{lemma}

\begin{proof}
In this case $P \cap P'$ is a single curve that meets $\Sigma_g$ in two points, and $D_i \cap D_i'$ is a single arc $\delta_i$ cutting out a subdisk $\Delta_i \subset D_i$ such that $\Delta_1 \cap \Sigma_g = \Delta_2 \cap \Sigma_g$.  Surgery on $D_i'$ along $\Delta_i$ yields a new compressing disk $D_i^+$, which can be chosen so that $\pd D_1^+ = \pd D_2^+$.  In addition $\pd D_i^+ \cap c = \pd D_i^+ \cap c' = \emp$.  Thus, $c^+ = \pd D_i^+$ is a reducing curve disjoint from both $c$ and $c'$.
\end{proof}

\begin{figure}[h!]
	\centering
	\includegraphics[width=.3\textwidth]{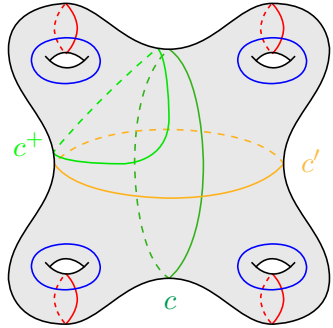}
\put (-67,10) {\textcolor{ForestGreen}{$c$}}
\put (-15,62) {\textcolor{Dandelion}{$c'$}}
\put (-124,68) {\textcolor{green}{$c^+$}}
	\caption{The case in which $\iota(c,c') = 2$}

	\label{2ps}
\end{figure}

Figure~\ref{2ps} depicts the curves from Lemma~\ref{2pts}.

A similar argument holds when $P \cap P'$ is two curves:

\begin{lemma}\label{2cvs}
If $\iota(c \cap c') = 4$ and $P \cap P'$ contains two curves, then $d_{\Hs}(c,c') = 2$.
\end{lemma}

\begin{proof}
In this case, each curve of $P \cap P'$ meets $\Sigma_g$ in precisely two points.  Thus, as in the proof of Lemma~\ref{2pts}, for $i=1,2$ there is an arc $\delta_i \subset D_i \cap D_i'$ cutting out a subdisk $\Delta_i \subset D_i$ with interior disjoint from $D_i'$, and such that $\Delta_1 \cap \Sigma_g = \Delta_2 \cap \Sigma_g$.  Let $D_i^{\pm}$ be the two disks obtained by surgery on $D_i'$ along $\Delta_i$.  Since $|\pd D_i^+ \cap \pd D| + |\pd D_i^- \cap \pd D| < 4$, we may suppose without loss of generality that $\pd D_i^+ \cap c = \emp$.  As above, we have $\pd D_1^+ = \pd D_2^+$, $\pd D_i^+ \cap c = \pd D_i^+ \cap c' = \emp$, and $c^+ = \pd D_i^+$ is a reducing curve disjoint from both $c$ and $c'$.
\end{proof}

In order to prove the next proposition, we first establish some facts about curves in a sphere with four boundary components, which we denote $\Sigma_0^4$.  Curves in $\Sigma_0^4$ are naturally parameterized by the extended rational numbers $\Q \cup \{\infty\}$, and given a rational number $a/b$, we let $\lambda_{a/b}$ denote the corresponding curve, where $\iota(\lambda_{a/b},\lambda_{c/d}) = 2|ad-bc|$.  The proof of the following lemma is elementary; see, for instance, Subsection 4.3 of~\cite{MZ} for further details.

\begin{figure}[h!]
	\centering
	\includegraphics[width=.3\textwidth]{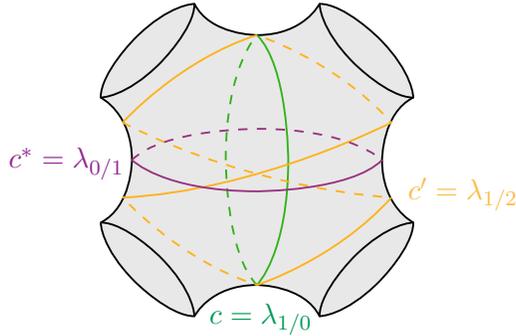}
\put (-83,3) {\textcolor{ForestGreen}{$c = \lambda_{1/0}$}}
\put (-8,50) {\textcolor{Dandelion}{$c' = \lambda_{1/2}$}}
\put (-159,62) {\textcolor{Plum}{$c^* = \lambda_{0/1}$}}
	\caption{Curves in $\Sigma_0^4$, where $\iota(c,c') = 4$ and $\tau_{c^*}(c) = c'$}

	\label{pillow}
\end{figure}

\begin{lemma}\label{param}
Suppose $c$ and $c'$ are two curves in $\Sigma_0^4$ such that $\iota(c,c') = 4$.  Then there is a parameterization such that $c = \lambda_{1/0}$ and $c' = \lambda_{\pm 1/2}$.  Additionally, if $c^* = \lambda_{0/1}$, then either $\tau_{c^*}(c) = c'$ or $\tau^{-1}_{c^*}(c) = c'$.
\end{lemma}
The curves $\lambda_{1/0}$, $\lambda_{0/1}$, and $\lambda_{1/2}$ are depicted in Figure~\ref{pillow}.

Next, we consider the remaining case, which is significantly more complicated than the previous two.

\begin{proposition}\label{4pts1cv}
If $\iota(c \cap c') = 4$ and $P \cap P'$ is a single curve, then $d_{\Hs}(c,c') < \infty$.
\end{proposition}

\begin{proof}
Observe that $D_1' \cap D_1$ is two arcs $\delta_1$ and $\delta_1^*$ that cobound disjoint subdisks $\Delta_1$ and $\Delta_1^*$ of $D_1$ with arcs $c_1$ and $c_1^*$ in $c$.  Similarly, $D_2' \cap D_2$ is two arcs $\delta_2$ and $\delta_2^*$ that cobound disjoint subdisks $\Delta_2$ and $\Delta_2^*$ of $D_2$ with arcs $c_2$ and $c_2^*$ in $c$.  Since $P \cap P'$ is a single curve, we have that $c = c_1 \cup c_2 \cup c_1^* \cup c_2^*$, where arcs meet only at their endpoints.  The setup is shown in Figure~\ref{forint}.  Let $E_1$ and $E_1^*$ be the result of surgery on $D_1'$ along $\Delta_1$.  We may assume that $E_1 \cap D_1 = \emp$ and $E_1^* \cap D_1 = \delta_1^*$.  We let $\Sigma^+$ and $\Sigma^-$ denote the two components of $\Sigma_g \setminus c$, and suppose without loss of generality that $\pd E_1 \subset \Sigma^+$.

\begin{figure}[h!]
	\centering
	\includegraphics[width=.3\textwidth]{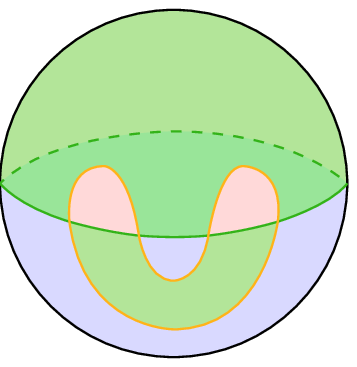}
\put (-134,64) {\textcolor{ForestGreen}{$c$}}
\put (-129,112) {$D_1$}
\put (-129,12) {$D_2$}
\put (-84,70) {\textcolor{Dandelion}{$\delta_1$}}
\put (-57,70) {\textcolor{Dandelion}{$\delta_1^*$}}
\put (-96,56) {\textcolor{red}{$\Delta_1$}}
\put (-47,56) {\textcolor{red}{$\Delta_1^*$}}
\put (-95,42) {\textcolor{ForestGreen}{$c_1$}}
\put (-47,40) {\textcolor{ForestGreen}{$c_1^*$}}
\put (-60,26) {\textcolor{Dandelion}{$\delta_2$}}
\put (-39,22) {\textcolor{Dandelion}{$\delta_2^*$}}
\put (-72,37) {\textcolor{blue}{$\Delta_2$}}
\put (-26,43) {\textcolor{blue}{$\Delta_2^*$}}
\put (-70,51) {\textcolor{ForestGreen}{$c_2$}}
\put (-70,90) {\textcolor{ForestGreen}{$c_2^*$}}
	\caption{Components of the intersections of $P'$ and $\Sigma_g$ with $P$}
	\label{forint}
\end{figure}

Now, surger $E_1^*$ along $\Delta_1^*$ to obtain disks $F_1$ and $G_1$, where $F_1 \cap D_1 = G_1 \cap D_1 = \emp$, and such that $\pd F_1 \subset \Sigma^+$ and $\pd G_1 \subset \Sigma^-$.  We also observe that the arcs $e_1$ and $f_1$ in $c' \cap \Sigma^+$ satisfy $\pd E_1 = e_1 \cup c_1$ and $\pd F_1 = f_1 \cup c_1^*$, and the arcs $g_1$ and $g_1^*$ in $c' \cap \Sigma^-$ satisfy $\pd G_1 = c_1 \cup g_1 \cup c_1^* \cup g_1^*$.  Here we are slightly abusing notation, since technically $\pd E_1$ is the union of $e_1$ and a slight pushoff of $c_1$ into $\Sigma^+$.  We note that both curves $\pd E_1$ and $\pd F_1$ are necessarily essential in $\Sigma_g$; otherwise, we could reduce $|c \cap c'|$ via isotopy.  However, it is possible that $\pd G_1$ is inessential.

We repeat a parallel construction in $H_2$:  The arcs $\delta_2$ and $\delta_2^*$ cobound disjoint subdisks $\Delta_2$ and $\Delta_2^*$ of $D_2$, and surgery on $D_2'$ along $\Delta_2$ and $\Delta_2^*$ yields disks $E_2$, $F_2$, and $G_2$ such that the curves $\pd E_2$ and $\pd F_2$ are essential, and $\pd E_2$ and $\pd F_2$ are contained in the same component of $\Sigma_g \setminus c$.  As above, we also suppose that $\pd E_2 = e_2 \cup c_2$ and $\pd F_2 = f_2 \cup c_2$, while $\pd G_2 = g_2 \cup c_2 \cup g_2^* \cup c_2^*$, where $e_2, f_2,g_2,g_2^*$ are the arcs of $c' \cap (\Sigma_g \setminus c)$.  It follows that the arcs $\{e_1,f_1,g_1,g_1\}$ are equal to $\{e_2,f_2,g_2,g_2^*\}$.  We already know the respective boundaries of each arc; thus, in pairs we have $\{e_1,f_1\} = \{g_2,g_2^*\}$ and $\{g_1,g_1^*\} = \{e_2,f_2\}$.  In other words, the disks $E_2$ and $F_2$ must be on the side of $c$ opposite $E_1$ and $F_1$, so that $\pd E_2 \cup \pd F_2 \subset \Sigma^-$ and $\pd G_2 \subset \Sigma^+$.  See Figure~\ref{setup}.

\begin{figure}[h!]
	\centering
	\includegraphics[width=.65\textwidth]{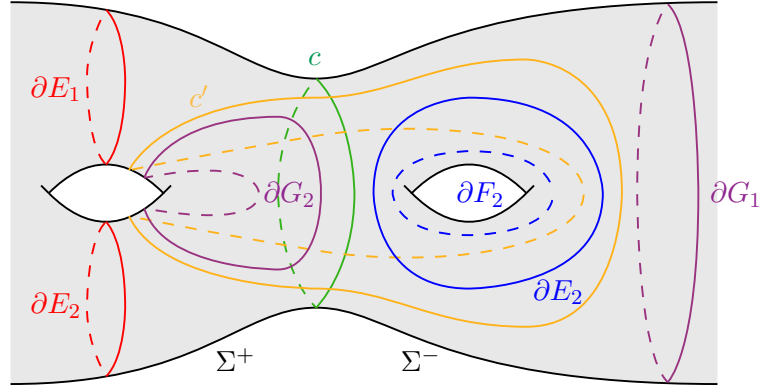}
\put (-162,128) {\textcolor{ForestGreen}{$c$}}
\put (-77,39) {\textcolor{blue}{$\pd E_2$}}
\put (-106,76) {\textcolor{blue}{$\pd F_2$}}
\put (-10,76) {\textcolor{Plum}{$\pd G_1$}}
\put (-179,76) {\textcolor{Plum}{$\pd G_2$}}
\put (-207,112) {\textcolor{Dandelion}{$c'$}}
\put (-267,116) {\textcolor{red}{$\pd E_1$}}
\put (-267,34) {\textcolor{red}{$\pd E_2$}}
\put (-197,12) {$\Sigma^+$}
\put (-127,12) {$\Sigma^-$}
	\caption{An example of boundaries of the disks obtained by surgeries on $D_1'$ and $D_2'$}
	\label{setup}
\end{figure}

Let $Q^+ \subset \Sigma^+$ be the subsurface $\overline{N}(c \cup e_1 \cup f_1)$.  Since $\pd e_1 = \pd c_1$ and $\pd f_1 = \pd c_1^*$, where $c_1$ and $c_1^*$ are disjoint arcs in $c$, it follows that the surface $Q^+$ is planar with four boundary components, one of which is the curve $c$.  The other three curves are parallel to $e_1 \cup c_1$, $f_1 \cup c_1^*$, and $e_1 \cup c_2 \cup f_1 \cup c_2^*$.  The first two curves are $\pd E_1$ and $\pd F_1$, respectively, and the third is $g_2 \cup c_2 \cup g_2^* \cup c_2^* = \pd G_2$.  Note that three of the boundary components of $Q^+$, namely $\pd E_1$, $\pd F_1$, and $c$, bound disks in $H_1$, so that the fourth boundary component, $\pd G_2$, also bounds a disk in $H_1$.  If $\pd G_2$ is an essential curve in $\Sigma_g$, then it bounds disks in both $H_1$ and $H_2$ and as such is a reducing curve in $\Sigma_g$ such that $c \cap \pd G_2 = c' \cap \pd G_2 = \emp$.  Thus, $d_{\Hs}(c,c') = 2$, completing the proof.  Otherwise, $\pd G_2$ bounds a disk in $\Sigma$, and we cap $Q^+$ off with this disk to obtain a pair of pants $R^+$ with boundary components $c$, $\pd E_1$, and $\pd F_1$.

We run a parallel construction in $\Sigma^-$:  Let $Q^- \subset \Sigma^-$ be the subsurface $\overline{N}(c \cup e_2 \cup f_2)$.  As above, $Q^-$ is a planar surface with four boundary components, $c$, $\pd E_2$, $\pd F_2$, and $\pd G_1$.  If $\pd G_1$ is essential in $\Sigma_g$, then it must be a reducing curve for $\Sigma_g$ and we have $d_{\Hs}(c,c') = 2$, as desired.  If not, then $\pd G_1$ bounds a disk in $\Sigma^-$, and we cap off $Q^-$ with the disk to obtain a pair of pants $R^-$ with boundary components $c$, $\pd E_2$, and $\pd F_2$.

If both $\pd G_1$ and $\pd G_2$ are inessential, let $R$ be the surface $R^+ \cup R^-$, so that $R$ is a sphere with four boundary components.  By construction, $R$ contains both curves $c$ and $c'$.  By Lemma~\ref{param}, we may choose a parameterization of curves in $R$ so that $c = \lambda_{1/0}$ and $c' = \lambda_{\pm 1/2}$.  Let $e \subset R$ be an arc from $\pd E_1$ to $\pd E_2$ that meets $c$ once and such that the lenses $\pd E_1$ and $\pd E_2$ and bridge $e$ determine an eyeglass twist $\psi$ with boundary $c^* = \lambda_{0/1}$.  By Lemma~\ref{eyfact}, we have
\[ \psi^{\pm 1}(c) = (\tau^{\pm 1}_{\pd E_1} \circ \tau^{\pm 1}_{\pd E_2} \circ \tau^{\mp 1}_{c^*})(c) = \tau^{\mp 1}_{c^*}(c),\]
so that by Lemma~\ref{param} either $\psi$ or $\psi^{-1}$ sends $c$ to $c'$.

If necessary, we may reverse the roles of $c$ and $c'$ to assume without loss of generality that $\psi(c) = c'$.  Consider an automorphism $\varphi$ of $\Sigma_g$ that sends $c$ to one of the curves $c_i^*$.  As in the proof of Lemma~\ref{eyenorm}, if $\psi_0$ is an eyeglass twist with lenses $\varphi(\pd E_1)$ and $\varphi(\pd E_2)$ and bridge $\varphi(e)$ that meets $c_i^* = \varphi(c)$ in a single point, then we have
\[ \varphi \circ \psi = \psi_0 \circ \varphi.\]
In addition, by Lemma~\ref{onebridge}, $\psi_0$ is a Powell move.  By Lemma~\ref{pact}, we have that $\psi_0(\Hs_0) = \Hs_0$, implying that $\varphi(c') = \varphi(\psi(c)) = \psi_0(\varphi(c)) = \psi_0(c_i^*) \in \Hs_0$, as is $\varphi(c) = c_i^*$.  It follows that $c'$ and $c$ are also in the same connected component, namely $\varphi^{-1}(\Hs_0)$, of $\Hs(\Sigma_g)$, completing the proof.
\end{proof}

\section{Reducing curves meeting in at most six points}\label{6p}

To extend our argument to reducing curves that meet in six points, and to prove Theorem~\ref{bigdist} in the following section, we employ a well-known tool, subsurface projection.  We say that a subsurface $\Sigma$ of the closed genus $g$ surface $\Sigma_g$ is \emph{essential} if $\Sigma$ is not an annulus or a pair of pants, and every boundary component of $\Sigma$ is essential in $\Sigma_g$.  Let $a$ be a properly embedded arc in $\Sigma$ such that $a$ is not isotopic to an arc in $\pd \Sigma$.  Then $P_a = \overline{N}(a \cup \pd \Sigma)$ is a pair of pants in $\Sigma$.  The \emph{subsurface projection} $\pi_{\Sigma}(a)$ is a subset of $\mathcal{C}(\Sigma)$ consisting of the curves of $\pd P_a$ that are essential in $\Sigma$.

Next, for any curve $c \subset \Sigma_g$, the \emph{subsurface projection} of $c$ to $\Sigma$ is defined to be the subset of $\mathcal{C}(\Sigma)$ given by the following conditions:
\be
\item If $c \subset \Sigma$, then $\pi_{\Sigma}(c) = \{c\}$.
\item If $c \cap \pd \Sigma \neq \emp$, then $\pi_{\Sigma}(c) = \bigcup_{a \subset c \cap \Sigma} \pi_{\Sigma}(a)$.
\item If $c \cap \Sigma =\emp$, then $\pi_{\Sigma}(c) = \emp$.
\ee
For further details, see~\cite{notes}.

As in Section~\ref{4p}, we set the convention that $c$ and $c'$ are reducing curves for $\Sigma_g$ such that $c$ bounds disks $D_1$ and $D_2$ in $H_1$ and $H_2$, respectively, and $c'$ bounds disks $D_1'$ and $D_2'$ in $H_1$ and $H_2$, respectively.  In addition, we let $P = D_1 \cup D_2$ and $P' = D_1' \cup D_2'$, isotoping $P$ and $P'$ to intersect minimally, and we let $\Sigma^{\pm}$ denote the two components of $\Sigma_g \setminus c$

\begin{lemma}\label{outermost}
Suppose that $c$ and $c'$ are reducing curves for $\Sigma_g$ such that $c \cap c' \neq \emp$.  Let $\delta$ be an arc of intersection of $D_i$ and $D_i'$ that is outermost in $D_i'$, where $\delta$ cobounds a subdisk $\Delta'$ of $D_i'$ with an arc $\alpha' \subset \pd D_i'$ with $\text{int}(\Delta') \cap D_i = \emp$.  Then either $\alpha' \subset \Sigma^+$ or $\alpha' \subset \Sigma^-$, and both curves in $\pi_{\Sigma^{\pm}}(\alpha')$ bound disks in $H_i$.
\end{lemma}

\begin{proof}
Since $\delta$ is outermost in $D_i'$, the arc $\alpha'$ meets $c$ only at its endpoints, so that $\alpha' \subset \Sigma^+$ or $\alpha' \subset \Sigma^-$.  Surgery on $D_i$ along $\Delta'$ yields disks $D_i^+$ and $D_i^-$, whose boundaries are precisely the two essential boundary curves of $\overline{N}(\alpha' \cup \pd \Sigma^{\pm})$ constituting $\pi_{\Sigma^{\pm}}(\alpha')$.
\end{proof}


Following Lemma~\ref{outermost}, we compare outermost arcs of intersection of $D_1 \cap D_1'$ and $D_2 \cap D_2'$.  For $i=1$ or $2$, let $\delta_i$ be an arc of $D_i \cap D_i'$ that is outermost in $D_i'$.  If $\delta_1$ and $\delta_2$ have the same endpoints, we say that the pair $P$ and $P'$ have \emph{matching bigons}.  If $\delta_1$ and $\delta_2$ have only one endpoint in common, we say that $P$ and $P'$ have \emph{adjacent bigons}.

\begin{lemma}\label{match}
Suppose $P$ and $P'$ have matching bigons.  Then there exists a reducing curve $c''$ such that $\iota(c, c'') =0$ and $\iota(c'',c') < \iota(c,c')$.
\end{lemma}

\begin{proof}
As in the proofs of Lemmas~\ref{2pts} and~\ref{2cvs} above, there are arcs $\delta_1 \subset D_1 \cap D_1'$ and $\delta_2 \subset D_2 \cap D_2'$ such that $\delta_i$ cuts out a subdisk $\Delta_i'$ of $D_i'$ whose interior misses $D_i$, and such that $\Delta_1' \cap \Sigma = \Delta_2' \cap \Sigma$.  In this case, surgery on $D_i$ along $\Delta_i'$ yields a compressing disk $D_i^+$ with that property that $D_i^+ \cap D_i = \emp$, $|\pd D_i^+ \cap \pd D_i'| < \iota(\pd D_i,\pd D_i')$, and $\pd D_1^+ = \pd D_2^+$.  Setting $c'' = \pd D_i^+$ completes the proof of the lemma.
\end{proof}

\begin{lemma}\label{adjac}
Suppose that $P$ and $P'$ have adjacent bigons.  Then there exists a reducing curve $c''$ such that $\iota(c,c'') = 4$ and $\iota(c'',c') < \iota(c,c')$.
\end{lemma}

\begin{figure}[h!]
\begin{subfigure}{.45\textwidth}
  \centering
  \includegraphics[width=.55\linewidth]{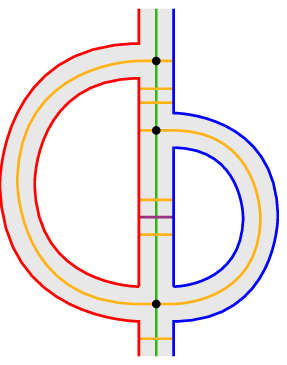}
\put (-51,16) {$x$}
\put (-51,128) {\textcolor{ForestGreen}{$c$}}
\put (-62,51) {\textcolor{Plum}{$e$}}
\put (-67,83) {$x_2$}
\put (-39,105) {$x_1$}
\put (-113,65) {\textcolor{Dandelion}{$\delta_1'$}}
\put (-4,50) {\textcolor{Dandelion}{$\delta_2'$}}

  \caption{Curve $c$, arc $e$, and arcs of $c'$ in $Q$}
  \label{adjA}
\end{subfigure}%
\begin{subfigure}{.45\textwidth}
  \centering
  \includegraphics[width=.55\linewidth]{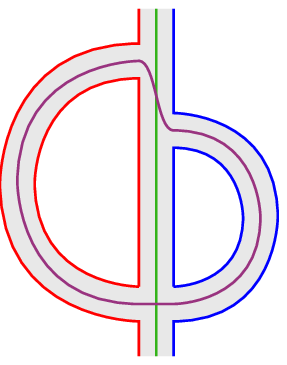}
  \put (-51,128) {\textcolor{ForestGreen}{$c$}}
  \put (-113,65) {\textcolor{Plum}{$c^*$}}
  \caption{Curves $c$ and $c^*$ in $Q$}
  \label{adjB}
\end{subfigure}\\
\begin{subfigure}{.45\textwidth}
  \centering
  \includegraphics[width=.55\linewidth]{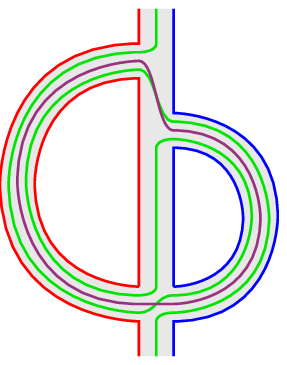}
    \put (-51,128) {\textcolor{green}{$c''$}}
  \put (-113,65) {\textcolor{Plum}{$c^*$}}
  \caption{Curves $c''$ and $c^*$ in $Q$}
  \label{adjC}
\end{subfigure}%
\begin{subfigure}{.45\textwidth}
  \centering
  \includegraphics[width=.55\linewidth]{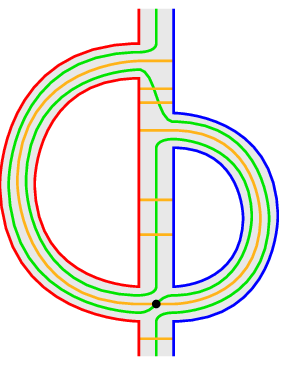}
      \put (-51,128) {\textcolor{green}{$c''$}}
      \put (-113,65) {\textcolor{Dandelion}{$\delta_1'$}}
\put (-4,50) {\textcolor{Dandelion}{$\delta_2'$}}
  \caption{Curve $c''$ and arcs of $c'$ in $Q$}
  \label{adjD}
\end{subfigure}\\
\caption{Four different pictures of curves $c$, $c^*$, and $c''$ along with arcs of $c'$ in the subsurface $Q$}
\label{adjaca}
\end{figure}

\begin{proof}
Since $\delta_1$ is outermost in $D_1'$, the arc $\delta_1$ cobounds a disk component of $D_1' \setminus D_1$ with an arc $\delta_1' \subset \pd D_1'$.  Similarly, $\delta_2$ cobounds a disk component of $D_2' \setminus D_1$ with an arc $\delta_2' \subset \pd D_2'$.  By assumption, $\delta_1'$ and $\delta_2'$ have one endpoint in common, call it $x$, which is contained in $c$, so we suppose without loss of generality that $\delta_1' \subset \Sigma^+$ and $\delta_2' \subset \Sigma^-$.  Let $x_1$ be the other endpoint of $\delta_1'$ and $x_2$ the other endpoint of $\delta_2'$.

Consider the subsurface $Q = \overline{N}(\delta_1' \cup \delta_2' \cup c)$, which is a sphere with four boundary components, $\pi_{\Sigma^+}(\delta_1')$ and $\pi_{\Sigma^-}(\delta_2')$, depicted in Figure~\ref{adjaca}.  By Lemma~\ref{outermost}, the two curves in $\pi_{\Sigma^+}(\delta_1')$ bound disks in $H_1$, while the two curves in $\pi_{\Sigma^-}(\delta_2')$ bound disks in $H_2$.  By construction, $c \subset Q$ and $c'$ meets $Q$ in the arc $\delta_1' \cup \delta_2'$, which intersects $c$ in the three points $x$, $x_1$, and $x_2$, and some number of additional arcs, each of which meets $c$ once, which we call \emph{short arcs} of $c' \cap Q$.  The setup is shown in Figure~\ref{adjA}.

Let $e$ be an arc in $Q$ that meets $c$ once and such that $e \cap c$ is contained in the arc component of $c \setminus (\delta_1' \cup \delta_2')$ with endpoints $x$ and $x_2$, as shown in Figure~\ref{adjA}.  Since $e$ connects a curve in $\D(H_1)$ to a curve in $\D(H_2)$, it determines an eyeglass twist $\psi$ with boundary curve $c^*$ shown in Figure~\ref{adjB}.  (We remark that if $e$ is instead chosen to meet $c$ between $x$ and $x_1$, it determines an eyeglass twist $\psi'$ which has different lenses but the same boundary curve $c^*$, so that the effects of $\psi^{\pm 1}$ and $(\psi')^{\pm 1}$ on $c$ are identical).

By Lemma~\ref{eyfact}, $\psi(c) = \tau_{c^*}(c)$, and thus there an eyeglass twist with boundary $c^*$ (or its inverse) that sends $c$ to another curve $c'' \subset Q$, such that $c''$ meets each of the short arcs of $c'$ once, and in addition $c''$ meets the arcs $\delta_1 \cup \delta_2$ in the single point $x$, instead of the three points of $c \cap c'$.  The curve $c''$ is shown first in Figure~\ref{adjC}, and its intersections with $c'$ are shown in Figure~\ref{adjD}.  We conclude that $c'' \in \Hs(\Sigma_g)$, $\iota(c,c'') = 4$, and $\iota(c'',c') < \iota(c,c'')$, completing the proof.
\end{proof}

\begin{remark}
We note that the combinatorics of the arcs described in the previous lemma determine Figure~\ref{adjaca} up to taking a mirror image.  Thus, in any case the statements in Lemma~\ref{adjac} are true for either an eyeglass twist $\psi$ or its inverse $\psi^{-1}$.
\end{remark}

\begin{proof}[Proof of Theorem~\ref{6pts}]
Suppose that $c,c' \in \Hs(\Sigma_g)$ with $\iota(c,c') \leq 6$.  If $\iota(c,c') = 2$, then $d_{\Hs}(c,c') = 2$ by Lemma~\ref{2pts}.  If $\iota(c,c') = 4$, then $d_{\Hs}(c,c') < \infty$ by Lemma~\ref{2cvs} and Proposition~\ref{4pts1cv}.  Thus, suppose that $\iota(c,c') = 6$.  As above, let $c$ bound disks $D_i$ in $H_i$, let $c'$ bound disks $D_i'$ in $H_i'$, and let $P = D_1 \cup D_2$ and $P' = D_1' \cup D_2'$ be the associated reducing spheres.  Note that $D_1 \cap D_1'$ has at least two arcs of intersection that are outermost in $D_1'$; pick two and call them $\delta_1$ and $\delta_1^*$.  Similarly, $D_2 \cap D_2'$ has at least two arcs of intersection $\delta_2$ and $\delta_2^*$ that are outermost in $D_2'$.  Since $\iota(c,c') = 6$ and the endpoints of the arcs $\delta_i$ and $\delta_i^*$ meet in points of $c \cap c'$, there must be a pair of arcs, say $\delta_1$ and $\delta_2$, with at least one common endpoint.

If $\delta_1$ and $\delta_2$ have both endpoints in common, then $P$ and $P'$ have matching bigons, and by Lemma~\ref{match}, there is a reducing curve $c''$ such that $\iota(c,c'') = 0$ and $\iota(c'',c') \leq 4$.  Otherwise, $P$ and $P'$ have adjacent bigons, and by Lemma~\ref{adjac}, there is a reducing curve $c''$ such that $\iota(c,c'') = 4$ and $\iota(c'',c') \leq 4$.  In either case, by the arguments mentioned above, $d_{\Hs}(c,c'') < \infty$ and $d_{\Hs}(c'',c') < \infty$, so that $d_{\Hs}(c,c') <\infty$, as desired.
\end{proof}

\begin{remark}
Note that the proof of the existence of matching or adjacent bigons fails when we consider $\iota(c,c') = 8$.  Indeed, it straightforward to construct reducing curves $c$ and $c'$ such that the corresponding reducing spheres $P$ and $P'$ do not have matching or adjacent bigons; thus, an inductive approach to proving that $\Hs(\Sigma_g)$ appears to fall short using only the methods given here.
\end{remark}

\section{Small intersection number but large distance in $\Hs(\Sigma_g)$}\label{large}

In this section, we prove Theorem~\ref{bigdist}, which asserts that for genus $g \geq 3$, reducing curves that meet in four points can be arbitrarily far apart in $\Hs(\Sigma_g)$, in a departure from the relationship between intersection number and distance in $\C(\Sigma_g)$ and $\D(H_i)$.  As a consequence we obtain Corollary~\ref{qi} about the geometry of $\Hs(\Sigma_g)$ as a subcomplex of $\D(H_i)$ and $\C(\Sigma_g)$.  Recall the definition of subsurface projection from the previous section.  The following lemma is well-known; see~\cite[Lemma 2.28]{notes}.

\begin{lemma}\label{projdist}
Let $\Sigma$ be an essential subsurface of $\Sigma_g$, and suppose that $c_0,c_1,\dots,c_k$ is a path in $\C(\Sigma_g)$ such that $c_i \cap \Sigma \neq \emp$ for all $i$.  Let $a \in \pi_{\Sigma}(c_0)$ and $a' \in \pi_{\Sigma}(c_k)$.  Then $d_{\C(\Sigma)}(a,a') \leq 6k$.
\end{lemma}

Suppose that $\Sigma$ is a surface with non-empty boundary, and let $\wh \Sigma$ denote the surface obtained by capping off each boundary component with a disk.   Note that any curve $c \subset \Sigma$ has a natural interpretation as a curve in $\wh \Sigma$, since $\Sigma$ includes into $\wh \Sigma$.  The next well-known lemma states that distance does not increase under this inclusion.

\begin{lemma}\label{closedist}
For any two curves $c,c' \in \C(\Sigma)$,
\[ d_{\C(\Sigma)}(c,c') \geq d_{\C(\wh\Sigma)}(c,c').\]
\end{lemma}

Recall that $\D(H_i)$ denotes the disk complex of $H_i$; that is, the subcomplex of $\C(\Sigma_g)$ induced by those curves that bound compressing disks in $H_i$.  For a curve $c \in \C(\Sigma_g)$, the distance from $c$ to $\D(H_i)$, denoted $d_{\C(\Sigma_g)}(c,\D(H_i))$, is
\[ d_{\C(\Sigma_g)}(c,\D(H_i)) = \min_{c' \in \D(H_i)} d_{\C(\Sigma_g)}(c,c').\]
The following theorem appears in work of Campisi and Rathbun~\cite[Theorem 1.2]{camprath}; another proof is based on work of Schleimer~\cite{saul}.
\begin{theorem}\label{fardist}\cite{camprath,saul}
Given a Heegaard splitting $S^3 = H_1 \cup_{\Sigma_g} H_2$ such that $g \geq 2$, and given $k \in \N$, there exists a curve $a \in \D(H_1)$ such that $d_{\C(\Sigma_g)}(a,\D(H_2)) \geq k$.
\end{theorem}

Note that if $S^3 = H_1 \cup_{\Sigma_g} H_2$ is a Heegaard splitting with a reducing curve $c$ that cuts $\Sigma_g$ into subsurfaces $\Sigma^+$ and $\Sigma^-$ with $g(\Sigma^{\pm}) = g^{\pm}$, then this decomposition induces genus $g^+$ and $g^-$ Heegaard splittings of $S^3$ in the following way:  Suppose that $c$ bounds disks $D_1 \subset H_1$ and $D_2 \subset H_2$.  Then the two components of $(H_1 \setminus D_1) \cup_{\Sigma_g \setminus c} (H_2 \setminus D_2)$ can be capped off with 3-balls, so that $\Sigma^{\pm}$ is capped off with a disk $D^{\pm}$, yielding Heegaard splittings we denote $H_1^+ \cup_{\wh\Sigma^+} H_2^+$ and $H_1^- \cup_{\wh\Sigma^-} H_2^-$.  We call these the \emph{Heegaard splittings induced by $c$}.  For any curve $c'$ in $\wh\Sigma^{\pm}$, we can choose a representative of $c'$ disjoint from $D^{\pm}$ and interpret $c'$ as a curve in $\Sigma^{\pm}$.  In particular, if $c'$ bounds a compressing disk in $H_i^{\pm}$, then it also bounds a disk in the original handlebody $H_i$.

\begin{lemma}\label{eyegl}
Suppose $S^3 = H_1 \cup_{\Sigma_g} H_2$ is a Heegaard splitting with reducing curve $c$ cutting $\Sigma_g$ into subsurfaces $\Sigma^{\pm}$, and let $\psi$ be an eyeglass twist with lenses $a \subset \Sigma^+$ and $b \subset \Sigma^-$ and bridge $e$ such that $|e \cap c| = 1$.  Then $c'=\psi(c)$ is a reducing curve such that $\iota(c,c') = 4$, $a \in \pi_{\Sigma^+}(c')$, and $b \in \pi_{\Sigma^-}(c')$.
\end{lemma}

\begin{proof}
This setup is shown in Figures~\ref{adjB} and~\ref{adjC} (with curve $c''$ in Figure~\ref{adjC} playing the role of $c'$ in this lemma).  By inspection, we verify that the claims of the lemma are true.
\end{proof} 

We have all the pieces in the place to prove the main theorem of this section.

\begin{proof}[Proof of Theorem~\ref{bigdist}]
Choose a reducing curve $c$ that cuts $\Sigma_g$ into subsurfaces $\Sigma^+$ and $\Sigma^-$ such that $g(\Sigma^+) = g-1$ and $g(\Sigma^-) = 1$.  Let $H_1^+ \cup_{\wh\Sigma^+} H_2^+$ and $H_1^- \cup_{\wh\Sigma^-} H_2^-$ be the Heegaard splittings induced by $c$.  By Theorem~\ref{fardist}, there exists a curve $a \in \D(H_1^+)$ such that $d_{\C(\wh\Sigma^+)}(a,\D(H_2^+)) \geq 6n$, and we fix a representative $a \subset \Sigma^+ \subset \wh \Sigma^+$, noting that $a \in \D(H_1)$.  Let $b$ be the unique curve in $\Sigma^-$ such that $b$ bounds a disk the solid torus $H_2^-$.  Let $e$ be an arc connecting $a$ to $b$ such that $|e \cap c|=1$, and let $\psi \in \G_g$ be the eyeglass twist with lenses $a$ and $b$ and bridge $e$.

Letting $c' = \psi(c)$, we have that $c'$ is another reducing curve for $\Sigma_g$, and since $|(a \cup b \cup e) \cap c| = 1$, it follows from Lemma~\ref{eyegl} that $\iota(c \cap c') = 4$.  In addition, since both $c$ and $c'$ are disjoint from the curves $a$ and $b$, we have
\[ d_{\D(H_i)}(c,c') = d_{\C(\Sigma_g)}(c,c') = 2.\]

For the final claim, suppose that $c' = c_0,c_1,\dots,c_m = c$ is a geodesic from $c'$ to $c$ in $\Hs(\Sigma_g)$, so that $c_{m-1} \cap c =\emp$ but $c_i \cap c \neq \emp$ for all $i < m-1$.  Additionally, the genus one surface $\Sigma^-$ does not contain a reducing curve disjoint from $c$; hence, $c_{m-1} \subset \Sigma^+$.  It follows that $c_i$ meets $\Sigma^+$ for all $i < m$.  By Lemma~\ref{eyegl}, we have that $a \in \pi_{\Sigma^+}(c_0)$, and $c_{m-1} \subset \Sigma^+$ implies $\pi_{\Sigma^+}(c_{m-1}) = \{c_{m-1}\}$.  Moreover $c_0,\dots,c_{m-1}$ is a path in $\C(\Sigma_g)$ such that every curve meets $\Sigma^+$, so by Lemma~\ref{projdist}, $d_{\Sigma^+}(a,c_{m-1}) \leq 6(m-1)$.  By Lemma~\ref{closedist}, we have $d_{\C(\wh\Sigma^+)}(a,c_{m-1}) \leq d_{\C(\Sigma^+)}(a,c_{m-1})$.  Finally, since $c_{m-1}$ is a reducing curve, it bounds a disk in $H_2^+$, which implies that $d_{\C(\wh\Sigma^+)}(a,\D_2(\wh\Sigma^+)) \leq d_{\C(\wh\Sigma^+)}(a,c_{m-1})$.  Combining these inequalities, we have
\[ 6n \leq d_{\C(\wh\Sigma^+)}(a,\D_2(\wh\Sigma^+)) \leq d_{\C(\wh\Sigma^+)}(a,c_{m-1}) \leq d_{\C(\Sigma^+)}(a,c_{m-1}) < 6m.\]
It follows that $n < m$, and since the path $c=c_0,\dots,c_m = c'$ in $\Hs(\Sigma_g)$ is a geodesic, we conclude that $n < d_{\Hs}(c,c')$, as desired.
\end{proof}

We need one final definition to prove the remaining corollary:  Given metric spaces $X$ and $Y$, a function $f:X \rightarrow Y$ is a \emph{quasi-isometric embedding} if there exist constants $A \geq 1$ and $B \geq 0$ such that for all $x_1,x_2 \in X$, we have
\[ \frac{1}{A} \cdot d_X(x_1,x_2) - B \leq d_Y(f(x_1),f(x_2)) \leq A \cdot d_X(x_1,x_2) + B.\]
It is well-known, for example, that the natural inclusion $\D(H_i) \hookrightarrow \C(\Sigma_g)$ is \emph{not} a quasi-isometric embedding (see Claim 4.12 of~\cite{notes} for a proof).  Corollary~\ref{qi} follows immediately from the combination of statements (2) and (3) of Theorem~\ref{bigdist}.

Despite the fact that the inclusion of $\D(H_i)$ into $\C(\Sigma_g)$ is not a quasi-isometric embedding, it is true that $\D(H_i)$ is \emph{quasi-convex} in $\C(\Sigma_g)$~\cite{MM2}.  In addition, the spaces $\C(\Sigma_g)$ and $\D(H_i)$ are known to be Gromov hyperbolic~\cite{MM,massch}.  Work of Akbas~\cite{akbas} implies that $\Hs(\Sigma_2)$ is quasi-isometric to a tree, so that $\Hs(\Sigma_2)$ is Gromov hyperbolic.  This leads us to two natural questions about the geometry of $\Hs(\Sigma_g)$:
\begin{question}
Is $\Hs(\Sigma_g)$ quasi-convex in $\D(H_i)$ or $\C(\Sigma_g)$?  Is $\Hs(\Sigma_g)$ Gromov hyperbolic?
\end{question}
Although these questions are most interesting in the event that $\Hs(\Sigma_g)$ is connected, recall that Corollary~\ref{3con} asserts that $\Hs(\Sigma_3)$ is connected, and it is our opinion that the Powell Conjecture is likely to be true, which would imply $\Hs(\Sigma_g)$ is connected for all $g$, lending merit to the questions above.

\bibliographystyle{amsalpha}
\bibliography{powbib}

\end{document}